\input amstex
\documentstyle{amsppt}
\magnification1200
\tolerance=10000
\def\n#1{\Bbb #1}
\def\p{\Bbb C_{\infty}}

\def\Gal{\hbox{Gal }}

\def\Hom{\hbox{Hom}}

\def\e11{E_{11}}
\def\gp{\goth p}
\def\ga{\goth A}

\def\vf{\varphi}

\def\vf{\varphi}

\def\de{\delta}

\def\ga{\gamma}
\def\la{\lambda}

\def\be{\beta}

\def\g{\goth }
\def\al{\alpha}

\def\th{\theta}
\def\vk{\varkappa}

\topmatter
\title SYSTEMS OF AFFINE EQUATIONS, THEIR L-FUNCTIONS, T-HOLONOMIC SEQUENCES, NON-COMMUTATIVE DETERMINANTS
\endtitle
\author
A. Grishkov, D. Logachev\footnotemark \footnotetext{E-mails: shuragri{\@}gmail.com; logachev94{\@}gmail.com (corresponding author)\phantom{****************}}
\endauthor
\thanks  The authors are grateful to V. Retakh who indicated them the
contents of Section 5.40 --- relations of non-commutative determinants
of the present paper and earlier notions of quasideterminants.
\endthanks
\NoRunningHeads
\address
First author: Instituto de Matem\'atica e estatistica, 
Universidade de S\~ao Paulo. Rua de Mat\~ao 1010, CEP 05508-090, S\~ao Paulo, Brasil, 
\medskip
Omsk State University n.a. F.M.Dostoevskii, Russia,
\medskip
and Sobolev Institute of Mathematics, Omsk, Russia.
\medskip
Second author: Departamento de Matem\'atica, Universidade Federal do Amazonas, Manaus, Brasil.
\endaddress
\abstract The authors defined in their previous paper
''$h^1\ne h_1$ for Anderson t-motives'' the notion of an
affine equation associated to a t-motive. We study in the present paper
general affine equations and their systems. We define L-functions of an
affine equation. We show that its local factor is an algebraic function;
we state a problem to prove an analog of the Lefschetz trace formula for
the global L-function. We define t-holonomic sequences which are
solutions of any affine equations. We state some conjectures and open
questions related to t-holonomic sequences, for example formulas of
their asymptotic growth. We show that some of the properties of
t-holonomic sequences are analogous to the properties of the classical
holonomic sequences, for example the closure properties. Further, we
show that systems of affine equations can be reduced to one equation by
a process analogous to the elimination of variables in a system of
linear equations; this gives us a notion of non-commutative
determinants. Finally, we consider some explicit examples of affine
equations that come from Drinfeld modules, their duals and some other
t-motives.
\endabstract
\endtopmatter

\NoRunningHeads

\document
{\bf 0. Introduction.} The Anderson ring $\p[T,\tau]$
of non-commutative polynomials in two variables $T$ and $\tau$ is
an analog of the ring of the linear differential operators, hence there
are analogies between the theory of Anderson t-motives (t-case) and the
theory of linear differential operators (D-case). First, there exist
t-case analogs of linear differential equations with polynomial
coefficients; we call them affine equations. Second, there exist t-case
analogs of classical holonomic sequences ([K]; [KP]); we call
them t-holonomic sequences.
\medskip
Affine equations were first defined in [GL21], Section 2, where we
started their study. Some affine equations come from Anderson t-motives.
\medskip
Here we continue to study affine equations and their systems. First, we
describe a process of non-commutative elimination of variables
(Proposition 5.8). As a result, we define a notion of non-commutative
determinant. By the way, the process of
non-commutative elimination of variables in a system of differential operators (Remark 5.37) apparently was not published
earlier; at least, it is not broadly known.
\medskip
Further, we formulate some open questions concerning t-holonomic
sequences, and we prove some of them in some particular cases. For
example, we have a question (here $v: \p^* \to \n Q^+$ is the valuation map on $\p$, and $x_0, \ x_1, \dots$ a t-holonomic sequence): 
\medskip
Is it true that all t-holonomic sequences are essentially of one of the
three types:
\medskip
A. $x_k\to 0$ (i.e. $v(x_k)\to +\infty$); 
\medskip
B. $x_k\to  c$ where $c$ is a constant, $c\in \n Q$ (maybe $c\in \n R$);
\medskip
C. $x_k\to -\infty$. 
\medskip
See Open question 3.11 for the exact statement.
\medskip
Further, we define a L-function of an affine equation (Section 1A). If
an affine equation comes from a t-motive, then it is the corresponding
L-function of this t-motive. We prove that its local factor is
''algebraic'' (Theorem 1A.1), and we ask to formulate and to prove an
analog of the Lefschetz trace formula (1A.1.5) for the global
L-function.
\medskip
Further, we appeal the reader to formulate and to prove an analog of the
Hilbert Nullstellensatz for affine equations (more exactly, for their
ideals) and for t-holonomic sequences, Questions 2.1. and 2.2.
\medskip
Let us recall that the Hilbert Nullstellensatz is the following. Let
$J\subset \p[x_1,\dots, x_n]$ be an ideal and
$V(J)\subset \n C^n$ the set of its zeros. For a subset
$Y\subset \n C^n$ let $I(Y)$ mean the set of polynomials
which are equal to 0 at $Y$ . Theorem: $I(V(J))$ is the radical of $J$. Hence, we want to know: What are
the relations between the ideal annihilating a t-holonomic sequence, and
the set of t-holonomic sequences annihilated by an ideal?
\medskip
Since the notion of a Gr\"obner basis exists for the Anderson ring
$\p[T,\tau]$, we can ask what are vertices of the left
monomial ideal attached to the ideal annihilating a t-holonomic
sequence. There is an analog of the same question for a system of affine
equations (0.3), (5.1).
\medskip
For the classical holonomic sequences we have a formula of its
asymptotical growth ([K], Theorem 2). We can expect that there
exists its t-case analog, i.e. a formula for the asymptotical value of
$v(x_k)$ (i.e. for $k>>0$) for t-holonomic
sequences of all the three above types. For
t-holonomic sequences of type (A) we 
give this conjectural formula in (3.7), we prove it for a particular
case in Propositions 3.2 (type (A)) and 3.8 (type (C)). An example of a
t-holonomic sequence of type B such that $v(x_k)=2-\frac3{2^{k+1}}$ is given in [GL21], Lemma 4.4a.
\medskip
t-holonomic sequences of type (A) are called small sequences. This is an
important notion: if we have an Anderson t-motive $M$ then its
$h^1(M)$ and $h_1(M)$ are equal to the dimension
of the space of small sequences of some two affine equations attached to
$M$, see [GL21].
\medskip
The classical holonomic sequences have some closure properties
([KP], Theorem 7.2). We can expect that these properties hold for
the t-case as well. We prove one of these properties for t-holonomic sequences (Theorem 4.5: the memberwise sum of two t-holonomic sequences is also t-holonomic).
\medskip
Also, we state two more research problems. First of them is given in
Section 4.4: is it possible to define a t-case analog of holonomic
functions of several variables ([K], Definition 3)?
\medskip
Second (Section 6.6): Let an affine equation $P$ of the rank
$r_0$ (or their system) be defined over a field $F\supset \n F_q(\th)$. In this case we have a representation $Gal(F)\to GL(X)$  where $X =\n F_q[[T]]^{r_0}$. What are the properties of this
representation?
\medskip
{\bf 0.1.} Let $M$ be an Anderson t-motive. We consider a system of affine equations defining $H^1(M)$, and (another) system of affine equations defining $H_1(M)$. The system defining $H_1(M)$ is the system defining $T_\g T(M)$ --- the $\g T$-Tate module of $M$, where $\g T$ is the prime ideal of $\n F_q[T]$ generated by $T$. See [NP] for a similar result. Further, the system defining $H^1(M)$ is the system defining $T_\g T(M^\vee)$, where $M^\vee$ is the 0-dual of $M$, see [HJ] for its definition. 
\medskip
Also, we show that the $T$-divisible modules of solutions of the
systems of affine equations defining first, $H_1(M)$, and second, $H^1(M)$, are distinct as $\Gal(\overline{\n F_q(\th)}/\n F_q(\th)$-modules (see (7.5): an example  for the case of Drinfeld modules of rank 2), although if $M$ is uniformizable then $H_1(M)$ and $H^1(M)$ are in canonical perfect duality.
\medskip
To understand a reason of this distinction is a subject of future research. Most likely we should find an equivalence relation on the set of systems of affine equations, as well as of the sets of these solutions, making them dual.
\medskip
{\bf 0.1.1.} Since the dimension of the set $T_\g T(M)_s$ --- the set of small elements (= sequences that tend to 0, see the definition below (2.2)) of $T_\g T(M)$, is $h_1(M)$, it is natural to ask what is the dimension of the set of small elements of $T_\g p(M)$ for all prime ideals $\g p$ of $\n F_q[T]$. The answer is simple: it is the same $h_1(M)$ for all $\g p$.
\medskip
Finally, there exists a natural embedding of $H_1(M)$ to $Lie(M)=\p^n$ where $n$ is the dimension of $M$. It defines a Hodge-Pink structure on $M$. Since $H_1(M)=T_\g T(M)_s$, we have an embedding $T_\g T(M)_s\hookrightarrow Lie(M)$.
\medskip
{\bf 0.1.2.} Do exist analogous embeddings of $T_\g p(M)_s$ for other prime ideals $\gp$?
\medskip
Generally, the present paper contains more questions than answers. 
\medskip
{\bf 0.2.} We use the standard definitions and notations for Anderson t-motives, see, for example, [A86], [G96], [GL20]. Namely, let $q$ be a power of a prime $p$, $\th$ a transcendental element, $\n R_\infty:=\n F_q((\th^{-1}))$, $\p:=\widehat{\overline{\n R_\infty}}$ an algebraically closed complete field with a valuation $v$ given by $v(\th)=-1$. The Anderson ring $\p[T,\tau]:=\p[T]\{\tau\}$ is a ring of non-commutative (i.e. $\tau a =a^q\tau$ for $a\in \p$) polynomials in two variables $T$ and $\tau$, and an Anderson t-motive $M$ is a left module over $\p[T,\tau]$ satisfying some properties. 
\medskip
We shall consider systems of equations $$\g P\cdot \g X=0\eqno{(0.3)}$$ (see (5.1), (5.2)) where $\g P\in M_{\la\times\mu}(\p[T,\tau])$ and $\g X\in M_{\mu\times1}(\p[[T]])$ a column of power series. We shall develop a theory of these systems generalizing the corresponding theory of [O33.1], [O33.2] and of [G96], Section 1 for the case $\la=\mu=1$ and entries of $\g P$, resp. of $\g X$ belong to $\p\{\tau\}$, resp. $\p$. Particularly, we describe an elimination of variables reducing the system (0.3) to the case $\la=\mu=1$, and we find the rank of the set of solutions for $\la=\mu=n$. 
\medskip
Some questions remain unsolved. Let $\g W$ be the set of solutions of (0.3)
and $\g W_1$ the set of solutions of a reduced system. What is an exact condition
that the projection from $\g W$ to $\g W_1$ is an isomorphism? What is a condition that two
systems of type (0.3) have a non-zero intersection of the sets of their
solutions (to treat both cases  $\la=\mu=1$ and $\la, \ \mu$ arbitrary)? Recall that for the case treated in [G86], [O33] this is a condition that the $p$-resultant is $\ne0$. See [E08] for a formula of $p$-resultant similar to the formula for ordinary resultant. 
\medskip
{\bf 1. t-holonomic sequences.} Let us consider the ring $(\p[[T]])\{\tau\}$ where the multiplication by $\tau$ is defined by the standard formula 
$$\tau(\sum_{i=0}^\infty a_iT^i) :=[(\sum_{i=0}^\infty a_iT^i)^{(1)}]\tau :=  [\sum_{i=0}^\infty a_i^qT^i]\tau\eqno{(1.1)}$$
\medskip
{\bf Proposition 1.2.} Its subring $\p[T,\tau]$ is left
noetherian.
\medskip
This is well-kniwn; for a convenience of the reader, let us give a
proof.
\medskip
{\bf Proof.} More generally, let $A$ be a commutative euclidean
ring with an euclidean degree function $d: A -\{0\} \to \n Z_{\ge 0}, d(0)=-\infty$, having an endomorphism $a \mapsto a^{(1)}$ 
such that $d(a)=d(a^{(1)})$. Let $A\{\tau\}$
be the ring of skew polynomials: $\tau a = a^{(1)}\tau$.
\medskip
Let us prove that $A\{\tau\}$ is left noetherian. Really, let $I \subset A\{\tau\}$ be a left ideal. Let $k$ be the minimal possible euclidean
degree of the leading coefficients of all
elements of $I$, and let $P=a_n\tau^n +...+a_0\in I$ be an element such that $d(a_n)=k$. 
\medskip
Applying the division algorithm we get immediately that for all
$m\ge n$, for all $X=b_m\tau^m +...+b_0\in I$ we have:
$a_n^{(m-n)}$ divides $b_m$. Really, let $b_m = a_n^{(m-n)}q+r$, where $d(r)<k$. We have 

$$X-q\tau^{m-n}P = r\tau^m+ ... + c_0\in I$$

hence $r=0$.
\medskip
Let us denote by $A\{\tau\}_{<n}$ the set of elements of $A\{\tau\}$ having the degree in $\tau$ less than $n$. Applying the above division construction to $X-q\tau^{m-n}P$ we get that $\forall \ X\in I \ \  \ \exists \ Q\in A\{\tau\}$, \ $R\in A\{\tau\}_{<n}$ such that
$$X=QP+R.$$

$A\{\tau\}_{<n}$ is a left $A$-module isomorphic
to $A^n$, hence $A\{\tau\}_{<n}\cap I$ is $A$-finitely generated. The
set of these generators and the element $P$ generate $I$ as a
left $A\{\tau\}$-module. $\square$
\medskip
The ring $\p[[T]]$ is a left module over $\p[T,\tau]$, where the multiplication by $\tau$ is defined by a version of the formula (1.1):
$$\tau(\sum_{i=0}^\infty a_iT^i) :=(\sum_{i=0}^\infty a_iT^i)^{(1)} :=  \sum_{i=0}^\infty a_i^qT^i\eqno{(1.3)}$$

Let $$X=\sum_{i=0}^\infty x_iT^i\in \p[[T]]\eqno{(1.4)}$$ be an unknown and $$P=\sum_{\ga=0}^{r_0} a_\ga\tau^\ga + \sum_{\be=1}^n\sum_{\ga=0}^{\vk_\be} b_{\be\ga} \tau^\ga T^\be\in\p[T,\tau]\eqno{(1.5)}$$ a coefficient. An affine equation is the equation $$PX=0,\eqno{(1.6)}$$see [GL21], (2.2.2). 
\medskip
Also, we can consider an analog of (1.6) for the case $P\in \p[[T]]\{\tau\}$, this is the equation [GL21] (2.1u), case of bounded $\vk_\be$. Or unbounded, to check.
\medskip
(1.6) is an infinite quantity of equations in $x_0, x_1, \dots$. For a convenience of the
reader, we give here the first of these equations, see also [GL21],
(2.1) - (2.1.2):
$$\hbox{The 0-th equation:} \ \ \ \ a_{r_0}x_0^{q^{r_0}}+a_{{r_0}-1}x_0^{q^{{r_0}-1}}+...+a_1x_0^q+a_0x_0=0\ \ \hbox{(head)}; \eqno{(1.6.0)}$$
(Here and below the set of terms $a_\ga x_i^{q^\ga}$ is called the head of the
equation, the set of other terms is called the tail of the equation. )
$$\hbox{The first equation:} \ \ \ \ a_{r_0}x_1^{q^{r_0}}+a_{{r_0}-1}x_1^{q^{{r_0}-1}}+...+a_1x_1^q+a_0x_1+\ \ \hbox{(head)} \eqno{(1.6.1)}$$ $$ +b_{1,\vk_1}x_0^{q^{\vk_1}}+ b_{1,\vk_1-1}x_0^{q^{\vk_1-1}}+...+b_{11}x_0^q+b_{10}x_0=0\ \ \hbox{(tail)};$$
$$\hbox{The second equation:} \ \ \ \ a_{r_0}x_2^{q^{r_0}}+a_{{r_0}-1}x_2^{q^{{r_0}-1}}+...+a_1x_2^q+a_0x_2+\ \ \hbox{(head)} \eqno{(1.6.2)}$$ $$ +b_{1,\vk_1}x_1^{q^{\vk_1}}+ b_{1,\vk_1-1}x_1^{q^{\vk_1-1}}+...+b_{11}x_1^q+b_{10}x_1+\ \ \hbox{(tail)}$$ $$ +b_{2,\vk_2}x_0^{q^{\vk_2}}+ b_{2,\vk_2-1}x_0^{q^{\vk_2-1}}+...+b_{21}x_0^q+b_{20}x_0=0;$$
etc., the length of the tail\footnotemark \footnotetext{The length of the tail is the quantity of unknowns in the tail part of the equation, i.e. it is 0, respectively 1, 2 in (1.6.0), respectively (1.6.1), (1.6.2) etc.} stabilizes, it does not exceed $n$.
\medskip
{\bf Definition 1.7.} A sequence $x_0, \ x_1, \dots$, where $x_i\in \p$, is called a t-holonomic sequence, if $X:=\sum_{i}x_iT^i$ is a solution of an above equation $PX=0$ ($P\in \p[T,\tau], \ P\ne0$).
By abuse of language,
$X$ as well is called a t-holonomic sequence.
\medskip
This is the t-case analog of the classical holonomic sequences (D-case;
we recall their definition in 4.1).
\medskip
 Let $P\in (\p[[T]])\{\tau\}$ be fixed. The set of $X\in \p[[T]]$ such that $PX = 0$ is a free $\n F_q[[T]]$-module of rank $r_0$ ([GL21],
Corollary 2.3.1).
\medskip
{\bf Remark 1.8.} The case $P\in (\p[[T]])\{\tau\}$ is not interesting,
because $\forall \ X \ \ \exists \ P\in (\p[[T]])\{\tau\}$ such that $PX=0$
([GL21], Proposition 2.3.2).
\medskip
{\bf 1.9.} $r_0$ is called the rank of $x_0, x_1, \dots$ and $n$ (the degree of $P$ in $T$) is
the length of its tail.
\medskip
{\bf Remark 1.9.1.} Let a t-holonomic sequence $X=\{x_0, x_1, \dots\}$ be given. The set of $P\in \p[T,\tau]$ such that $PX=0$ is a left ideal denoted by $I(X)$. Different generators of $I(X)$ can have different $r_0$ and $n$, hence the above (1.9) must be modified.
\medskip
{\bf Remark 1.10.} From the first sight, we can always transform
$P$ in some $P_0$ such that for $P_0$ we have: all $\vk_\be<r_0$. But after this modification we shall get: $P_0\in (\p[[T]])\{\tau\}$, i.e. $P_0\not\in\p[T,\tau]$. Example: let $$P=\tau+a_0+(\tau+b_{10})T;$$
The equations (1.6.*) are
$$x_0^q+a_0x_0=0; \ \ \ \matrix x_1^q+a_0x_1+ \\ x_0^q+b_{10}x_0=0; \endmatrix \ \ \ \matrix x_2^q+a_0x_2+ \\ x_1^q+b_{10}x_1=0;  \endmatrix \ \ \ \hbox{etc.}$$
The transformation will give us 
$$x_0^q+a_0x_0=0; \ \ \ \matrix x_1^q+a_0x_1+ \\ (b_{10}-a_0)x_0=0; \endmatrix \ \ \ \matrix x_2^q+a_0x_2+ \\ (b_{10}-a_0)x_1-\\ (b_{10}-a_0)x_0=0;  \endmatrix \ \ \ \matrix x_3^q+a_0x_3+ \\ (b_{10}-a_0)x_2-\\ (b_{10}-a_0)x_1+\\ (b_{10}-a_0)x_0=0;  \endmatrix \ \ \ \hbox{etc.}$$
i.e. 
\medskip
$P_0=\tau+a_0+(b_{10}-a_0)T-(b_{10}-a_0)T^2+(b_{10}-a_0)T^3-...\in (\p[[T]])\{\tau\}$. 
\medskip
We have $(1+T)P_0=P$, i.e. this modification
corresponds to a multiplication by an element of $(\p[[T]])\{\tau\}^*=\p[[T]]^*$.
\medskip
Non-t-holonomic sequences exist, i.e. this notion is non-trivial. The
simplest way to prove this fact is to show that $x_i$ cannot grow
too rapidly:
\medskip
{\bf Proposition 1.11.} If $\forall \ m$ we have $v(x_m)=-2^{m^2}$ then $x_0, x_1, \dots$ is not t-holonomic. 
\medskip
{\bf Proof.} We have: $a_{r_0}x_m^{q^{r_0}}$ is the leading term of
(1.6.0). It is sufficient to show that $\forall \ P \ \ \ \exists \ m$ such that $v(a_{r_0}x_m^{q^{r_0}})$ is
strictly less than $v$ of all other terms of (1.6.*), i.e. that 
$$\forall \ P \ \ \ \exists \ m \ \ \ \forall \ \be, \ \ga \hbox{  we have } v(a_{r_0}x_m^{q^{r_0}})<v(a_{\ga}x_m^{q^{\ga}}) \ (\ga\ne r_0);$$
$$v(a_{r_0}x_m^{q^{r_0}})<v(b_{\be\ga}x_{m-\be}^{q^{\ga}}).$$
This is immediate. $\square$
\medskip
Moreover, ''almost all'' sequences are non-t-holonomic. (Indicate what
''almost all'' mean).
\medskip
{\bf 1A. L-functions of affine equations.}
\medskip
There are several types of L-functions of Anderson t-motives. Let us
recall a definition of the type that we shall generalize.
\medskip
Let us fix a number $d\in \n Z_{\ge1}$. Let us consider first the local case, i.e. let $M$ be an Anderson t-motive of rank $r$ over $\n F_{q^d}$, i.e. a module
over $\n F_{q^d}[T,\tau]$ satisfying 
standard conditions. We consider $M$ as a free module of dimension
$r$ over $\n F_{q^d}[T]$, with a skew action of $\tau$.
\medskip
For a matrix $A=\{a_{ij}\} \in M_{r\times r}(\n F_{q^d}[T])$ we denote by $A^{(k)}$ the matrix with entries $\{a_{ij}^{(k)}\}$, and $A^{[k]}:=A^{(k)-1}\cdot ...\cdot A^{(1)}\cdot A$. 
\medskip
Let $f_*=(f_1, \dots, f_r)^t$ (a matrix column; $^t$ means transposition) be a basis of $M$ over $\n F_{q^d}[T]$. Let the action of $\tau$ on $f_*$ be given by a matrix $Q\in M_{r\times r}(\n F_{q^d}[T])$, i.e. $$\tau f_*=Qf_*.$$ We have: 
$$L(M,U):=\det (I_r-Q^{[d]}U^d)^{-1}\in \n F_{q}[T][[U]]$$
(because obviously $\det (I_r-Q^{[d]}U^d)\in \n F_{q}[T,U]$).
\medskip
The global L-function is a product of almost all Euler factors (we
neglect here the bad primes). Namely, let $K$ be a finite extension of $\n F_{q}(\th)$, \ $M$ an Anderson
t-motive of rank $r$ over $K$, $\g P$ a prime of $K$ such that
the reduction of $M$ at $\g P$ is
good. Let the residue field of $\g P$ be $\n F_{q^d}$. Let $L_\g P(M,U)$ be the local L-function of
the reduction of $M$ at $\g P$. The global L-function up to finitely
many Euler factors is
$$L(M,U):=\underset{\g P}\to{\prod}L_\g P(M,U).$$

Also, we can consider the characteristic polynomial of the action of
Frobenius
on the Tate module of $M$. Let us consider the local case, i.e.
$M$ is over $\n F_{q^d}[T]$. Let $\g L$ be a prime of
$\n F_{q}[T]$. The $\g L$-Tate module $T_\g L(M)$ of
$M$ is defined (where? is it on an analog of $E(M)$?). Frobenius (which one?) acts on $T_\g L(M)$, and its characteristic polynomial
(which one?) coinsides with the characteristic polynomial of $Q$ 
(why?)
\medskip
Let $P\in \n F_{q^d}[T,\tau]$. We consider the equation
(1.6) $PX=0$ for this $P$ (notations from (1.5)), where $X\in \overline{\n F_q}[[T]]$. As above, let $V(P)\subset \overline{\n F_q}[[T]]$ be the set of solutions
of $PX=0$, it is (as earlier) a free $\n F_q[[T]]$-module of dimension $r=r_0$.
\medskip
The Frobenius automorphism $Fr: \sum x_iT^i \mapsto  \sum x_i^{q^d}T^i$ is an
automorphism of $V(P)$ (we have $Fr=\tau^d$). Its
characteristic polynomial $Ch(P):=|Fr-Id_rU|$  belongs to $\n F_q[[T]][U]$. 
We denote its coefficients by $A_i$, i.e. $Ch(P)=\sum_{i=0}^rA_iU^i$. 
\medskip
{\bf Theorem 1A.1.} All $A_i$ are elements of $\n F_q[[T]]$ 
 which are algebraic over $\n F_q(T)$. 
\medskip
{\bf Proof.} Let us consider first the case $d=1$. Let $X_1, X_2, \dots,X_r$ be a
$\n F_q[[T]]$-basis of $V(P)$. We have $X_i\in \overline{\n F_q}[[T]]$. 
\medskip
Let us consider $P$ as a polynomial in $\tau$, i.e. $$P=\sum_{i=0}^\vk S_i\tau^i\eqno{(1A.1.1)}$$
where $S_i\in \n F_q[T]$ (we have $\vk = max(\vk_\be)$ from (1.5)). Let $\g X := (X_1,X_2, \dots, X_r)^t$ be this basis written as a matrix column. We denote
by $C\in M_{r\times r}(\n F_q[[T]])$ the matrix of the
action of $\tau = Fr$ on $\g X: \tau(\g X)=C\g X$. 
\medskip
We have $\tau^i(\g X)=C^i\g X$, hence using notations of (1A.1.1) we
have $$[\sum_{i=0}^\vk S_iC^i](\g X)=0\eqno{(1A.1.2)}$$ i.e. $\sum_{i=0}^\vk S_iC^i=0$. 
\medskip
Let $\la_1, \dots, \la_\ga$ be all different (i.e. $\la_i\ne \la_j$) eigenvalues of $Fr$ considered as a linear
transformation on $V(P)$. We have $\la_i\in \overline{\n F_q((T))}$. Let $m_i$ be the maximal size of a Jordan
block of $\la_i$, and $k_i$ the total size of all Jordan blocks of
$\la_i$, i.e. the minimal polynomial of $Fr$ 
is $\prod_{i=1}^\ga(U-\la_i)^{m_i}$ (up to a sign).
\medskip
We have $\prod_{i=1}^\ga(U-\la_i)^{m_i}$ divides $\sum_{i=0}^\vk S_iU^i$. Since all $S_i$ belong to $\n F_q[T]$, we get that
$\la_i\in \overline{\n F_q(T)}$. Since $Ch(P)=\prod_{i=1}^\ga(U-\la_i)^{k_i}$ (up to a sign), we get that $A_i\in \overline{\n F_q(T)}$.
\medskip
For the case of any $d$ we need a lemma. Let $K/K_0$ be a
cyclic Galois extension of fields, $Gal(K/K_0=\n Z/d\n Z$. We fix a generator $\ga$ of $Gal(K/K_0$. Let
$Z$ be an independent variable. Let $K\{Z\}$ be a skew field of
non-commutative polynomials $$P=\sum_{i=0}^n x_iZ^i, \ \ \ x_i\in K, \eqno{(1A1.3)}$$
where the multiplication is defined the formula (here $x\in K$):
$$Zx=\ga(x)Z.$$
The field $K_0[Z^d]$ is the center of $K\{Z\}$, the
dimension of $K\{Z\}$ over $K_0[Z^d]$ is $d^2$.
\medskip
{\bf Lemma 1A.1.4.} Let $P$ be from 1A.1.3. There exists the
only $Q\in K\{Z\}$ such that $QP\in K_0[Z^d]$ is of
degree $nd$, i.e. $QP=\sum_{i=0}^n y_iZ^{id}$, $y_i\in K_0$, and $y_n=N_{K/K_0}(x_n)$. 
\medskip
For a proof see any textbook on skew fields. The element $QP$ is
called the reduced norm of $P$. (Maybe $K$ must contain
$\zeta_d$ - a primitive root of 1 of order $d$?)
\medskip
We apply this lemma to the case $K=\n F_{q^d}(T), \ K_0=\n F_{q}(T), \ \ga(x)=x^{(1)}$ for $x\in \n F_{q^d}(T)$, and $Z=\tau$ in the skew field $\n F_{q^d}(T)\{\tau\}$. 
\medskip
For $P$ from (1A.1.1) we have $S_i\in K$, i.e. $P$ from (1A.1.1) is $P$ from (1A.1.3). Let $C$ be as above, but now $Fr=\tau^d$, i.e. $\tau^d(\g X)=C\g X$. 
\medskip
Let $Q, \ QP, \ y_i$ be from the lemma. We have $QP(\g X)=0$, i.e. $[\sum_{i=0}^ny_iC^i](\g X)=0$. 
This is an analog of (1A.1.2),
$y_i$ are $S_i$, and $n$ is $\vk$. The end of the proof for the
case of any $d$ coincides with the above proof for $d=1$. $\square$
\medskip
{\bf Remark.} The proof shows that if for all $i$ we have
$\vk_i<r$ then $A_I\in \n F_q[T]$. If for
all $i$ we have $\vk_i\le r$ then $A_I\in \n F_q(T)$.
Clearly the analog of this theorem exists for systems of affine
equations 5.1.
\medskip
Let us recall the Lefschetz trace formula for the L-functions of
Anderson t-motives. It is formulated in terms of $\Cal E$-$\tau$-sheaves. We have:

$$L(\Cal E, \tau, U)=\frac{\det(1-H^1(\n P^1, \tau\circ\vf)\cdot U)}{\det(1-H^0(\n P^1, \tau\circ\vf)\cdot U)}\eqno{(1A.1.5)}$$ This formula implies that $L(\Cal E, \tau, U)\in \n F_q[t](U)\cap  \n F_q[t][[U]]$. (To check notations: $t$, $T, \ U$). 
\medskip
{\bf Problem.} To formulate and to prove an analog of (1A.1.5) for
affine equations. Is it true that $L(P,U)\in \overline{\n F_q(t,U)}\cap  \n F_q[t][[U]]$?
\medskip
{\bf 2. Power series and ideals: analog of the Hilbert's
Nullstellensatz.}
\medskip
Here we state two questions which form an analog of the Hilbert's
Nullstellensatz for the t-case. The answers are subjects of further research.
\medskip
Let $X$ be a t-holonomic sequence. The set of $P\in \p[T,\tau]$ such that $PX=0$ is a left ideal of
$\p[T,\tau]$. We denote it by $I(X)$. More
generally, for a set of t-holonomic sequences $Z$ we denote by
$I(Z)$ the left ideal consisting of all $P\in \p[T,\tau]$ such that $\forall \ X \in Z$ we have $PX=0$. 
\medskip
Let $J\subset\p[T,\tau]$ be a left ideal. We denote by
$V(J)$ the set of $X\in\p[[T]]$
such that for all $Q\in J$ we have $QX=0$. Obviously $\n F_q[[T]]\cdot X\subset V(I(X))$. 
\medskip
$V(I(X))$ can be strictly greater than $\n F_q[[T]]\cdot X$, see a below example. 
\medskip
{\bf Question 2.2.} What is the relation between $J$ and
$I(V(J))$?
\medskip
Later (in this section) we shall consider $\p[T,\tau]$ as
the ring of polynomials in $\tau$ with coefficients in
$\p[T]$, and the word ''degree'' will mean the degree of
$P\in\p[T,\tau]$ as a polynomial in $\tau$.
\medskip
{\bf Lemma 2.3.} Let $J\subset\p[T,\tau]$ be a left ideal. There exists the only (up
to a factor from $\p^*$) polynomial $P_0\in J$ of the degree $d$ such
that:
\medskip
1. All $P\ne 0 \in J$ have the degree $\ge d$;
\medskip
2. All $P\in J$ of the degree $d$ have the form $P=QP_0$ where $Q\in \p[T]$. 
\medskip
{\bf Proof} is immediate. We choose $P_0$ as an element of $J$ having first, the minimal
degree (we denote it by $d$), and second, its coefficient at $\tau^d$ (it belongs to $\p[T]$)
has the minimal degree in $T$. Properties (1), (2) trivially hold. $\square$
\medskip
For a left ideal $J\subset\p[T,\tau]$  we denote the polynomial $P_0$ from Lemma 2.3 by $min(J)$. 
\medskip
{\bf Definition 2.4.} A polynomial $P\subset\p[T,\tau]$ is called non-reducible if its coefficients have no common
factor, i.e. if $Q\in\p[T]$ divides $P$ from the left (i.e. $\exists \ P'$ such that $P=QP'$) then
$Q$ is constant.
\medskip
{\bf Remark 2.8.} Let $J\subset\p[T,\tau]$ be a left ideal. $min(J)$ is not necessarily non-reducible. 
\medskip
{\bf Definition 2.4a.} A left ideal $J\subset\p[T,\tau]$ is called saturated if for any $P\in J$, for any $Q\in\p[T]$ we have: $\{\ Q^{-1}P\in \p[T,\tau] \ \} \implies \{\ Q^{-1}P\in J\ \}$. 
\medskip
{\bf Definition 2.5.} Let $J\subset\p[T,\tau]$ be a
left ideal. The set of all elements of the form $\{\ Q^{-1}P\ \}$ where $P\in J$, $Q\in\p[T]$ are such that $Q^{-1}P\in\p[T,\tau]$, is called the saturation of $J$. 
\medskip
Obviously the saturation is a left ideal. A left ideal is saturated iff it
coincides with its saturation. We denote the saturation of a left ideal $J$ by $\bar J$.
\medskip
{\bf Remark 2.5a.} It is possible that $P_0$ 
is non-reducible but $\p[T,\tau]\cdot P_0$ is not saturated. Example: 
$$P_0=(T+\th)\tau+(T+\th^q), \ \ \ \ P=\tau^2+\frac{\th^{q^2}-\th}{\th^{q}-\th}\tau+\frac{\th^{q^2}-\th^q}{\th^{q}-\th}.$$
We have: 
$$(\tau+\frac{\th^{q^2}-\th^q}{\th^{q}-\th})P_0=(T+\th^q)P\in \p[T,\tau]\cdot P_0\eqno{(2.5a.1)}$$
but $P\not\in \p[T,\tau]\cdot P_0$. This is a
consequence of non-commutativity: $\tau\cdot(T+\th)=(T+\th^q)\tau$.
\medskip
{\bf Lemma 2.6.} Let $J\subset\p[T,\tau]$ be a left
ideal. Then $J\subset\overline{\p[T,\tau]\cdot min(J)}$. 
\medskip
{\bf Proof.} Trivial. We must prove that if $P\in J$ then $\exists \ Q\in \p[T]$ such that $QP\in \p[T,\tau]\cdot P_0$. We use induction by the degree of $P$ in $\tau$. Let $P_0=C_d\tau^d+C_{d-1}\tau^{d-1}+...+C_0$ where $C_i\in \p[T]$. If the
degree of $P$ in $\tau$ is $<d$ then $P=0$. Let the above $Q$ exists for all $P$ whose degree in
$\tau$ is $\le m$, and let us prove its existence for $P$ whose
degree in $\tau$ is $m+1$. Let $P=K_{m+1}\tau^{m+1}+...$. We have 
$$C_d^{(m+1-d)}P-K_{m+1}\tau^{m+1-d}P_0\in J,$$
and its degree in $\tau$ is $\le m$, i.e. by induction $\exists \ Q_m\in \p[T]$ such that
$$Q_m(C_d^{(m+1-d)}P-K_{m+1}\tau^{m+1-d}P_0)=S_mP_0$$

for some $S_m\in \p[T,\tau]$. This means that
$$Q_mC_d^{(m+1-d)}P\in \p[T,\tau]\cdot P_0. \ \ \ \ \ \ \ \ \ \ \ \ \ \square$$

{\bf Corollary 2.7.} We have $\bar J=\overline{\p[T,\tau]\cdot P_0}$.
\medskip
{\bf Observation 2.9.} Let $J\subset \p[T,\tau]$ be a
left ideal. When $I(V(J))=\bar J$?
\medskip
Not alwais. Idea: Let $P\in I(V(J))$. Let us consider
the left ideal $<P,J>:=\{\p[T,\tau]\cdot P+J$, and let $P_0=min(J), \ P_1=min(<P,J>)$. If
$P_1=P_0$, then $P\in \overline{\p[T,\tau]\cdot P_0}$, i.e. $P\in \bar J$. 
\medskip
If $deg(P_1)=deg(P_0)$ then $P_0=UP_1$ for some $U\in \p[T,\tau]$. We have $\overline{\p[T,\tau]\cdot P_0}=\overline{\p[T,\tau]\cdot P_1}$. Really, let $P_1=\sum_{i=0}^dC_i\tau^i$ where $C_i\in \p[T]$, and let $X=(\sum_{i=0}^kB_i\tau^i)P_1\in\p[T,\tau]\cdot P_1$. We must prove that $\exists \ Q \in \p[T]$ such that $QX\in\p[T,\tau]\cdot P_0$. We have $Q=U\cdot U^{(1)}\cdot ... \cdot U^{(d)}$, because $QX=(\sum_{i=0}^kK_i\tau^i)P_0$ for $K_i=U\cdot U^{(1)}\cdot ... \widehat{U^{(i)}} ... \cdot U^{(d)}\cdot B_i$. 
\medskip
What to do if $deg(P_1)<deg(P_0)$? 
\medskip

Ukazat', chto takoe $n$ str. 5, pered Definition 1.7.
\medskip
{\bf Remark 2.9.0}: (see also 1.9): if $x_0, x_1,...$ are given, then $P$ may be different, i.e. $r_0$ and $n$ are not unique as well. 
\medskip
We denote $A:=\p[T,\tau]$. 
\medskip
{\bf Remark 2.9a.} Lemma 2.6 means that any left ideal of $A$ is "near to principal".
\medskip
{\bf Remark 2.9b.} Statement of Proposition 2.9 is false. Counterexample: $J=A\cdot \tau$. We have $V(J)=0$, $I(V(J))=A\ne\overline{A\cdot \tau}$.
\medskip
I think that Proposition 2.9 is "near to correct".
\medskip
{\bf Question 2.10.} Does exist $X$ such that $I(X)$ is not a principal left ideal?
\medskip
{\bf Answer.} Yes. Let $P_0$ be from 2.5a, and $X$ satisfies $P_0X=0$. (2.5a.1) implies that $PX=0$, i.e. $P_0, \ P\in I(X)$. It is clear that $I(X)$ is not a principal left ideal. Really, let exist a generator $Q$ of $I(X)$. Then $\exists \ R_0, \ R$ such that $P_0=R_0Q, \ P=RQ$. The degrees in both $T$, $\tau$ are additive in multiplication, hence $deg_T(Q)=0$, $deg_\tau(Q)=0$ or 1. Clearly $deg_\tau(Q)\ne0$. If $deg_\tau(Q)=1$, i.e. $Q=x\tau+y$, then a simple calculation shows that $P_0=R_0Q$ does not hold. $\square$
\medskip
{\bf Definition 2.11.} $X=\sum_{i=0}^\infty x_iT^i\in \p[[T]]$ is called bounded iff there exists a finite $\n F_q$-module $ G\subset \p$ such that $\forall\ i$ we have $ x_i\in G$. 
\medskip
Obviously bounded are holonomic. 
\medskip
If $X$ is bounded then $I(X)\cap \p\{\tau\}\ne0$. Conversely, if $J\cap \p\{\tau\}\ne0$ then all elements of $V(J)$ are bounded. 
\medskip
$X$ from Question 2.10 is bounded, because $P\in \p\{\tau\}$. 
\medskip
{\bf Preliminary calculation.} I think that $\exists \ P_0, \ P_1\in A$ of degrees: 
\medskip
$deg_T(P_0)=2$, $deg_\tau(P_0)=1$
\medskip
$deg_T(P_1)=1$, $deg_\tau(P_1)=2$
\medskip
such that $\exists \ R_0=\tau+x, \ R_1=T+w$ ($x, \ w\in \p$) such that 
$$R_0P_0=R_1P_1\eqno{(2.12)}$$ and such that the ideal $AP_0+AP_1$ is not principal. 
\medskip
Let $X\in \p[[T]]$ satisfies $P_0X=0$. (2.12) implies that $P_1X=0$. We can expect that $I(X)$ is not principal, and $X$ is not bounded. 
\medskip
{\bf Question 2.13.} Does exist $\g X\subset \p[[T]]$ such that $\g X$ is isomorphic to $(\n F_q[[T]])^r$, such that for any $P\in I(\g X)$ we have: The set of $Y\in \p[[T]]$ such that $PY=0$ is strictly greater than $\g X$? 
\medskip
Let us consider, for example, $P_1=\tau^2+\th + T$. If $X=\sum_{i=0}^\infty x_iT^i$ such that $P_1X=0$ then $$x_0^{q^2}+\th x_0=0$$ $$x_1^{q^2}+\th x_1+x_0=0\hbox{ etc.}$$Let $X_0$ be a fixed solution of $P_1X_0=0$. 
\medskip
{\bf Proposition 2.14. } $I(X_0)=\overline{AP_1}.$
\medskip
{\bf Proof.} Let $P_0=min(I(X_0))$. If $P_0=P_1$ then the proposition is clear (proveril, eto tak). Clearly $deg_\tau(P_0)\ne0$, hence we consider the case $deg_\tau(P_0)=1$. Let $P_0=Y_1\tau+Y_0$, $Y_i\in \p[T]$. According Lemma 2.6, $\exists \ Q\in \p[T]$ such that $$QP_1=(Z_1\tau+Z_0)P_0,\eqno{(2.15)}$$ $Z_i\in \p[T]$. (2.15) is $$Q=Z_1Y_1^{(1)}$$ $$Z_0Y_1+Z_1Y_0^{(1)}=0\eqno{(2.16)}$$ $$Q(T+\th)=Z_0Y_0$$Let $deg_T(Z_i)=z_i, \ deg_T(Y_i)=y_i$. Equations (2.16) imply $$z_1+y_1+1=z_0+y_0, \ \ \ z_0+y_1=z_1+y_0.$$These equations are incompatible modulo 2 --- a contradiction. $\square$
\medskip
We get that the answer to Question 2.13 is "yes". Really, we choose such $X_0$, and $\g X=\n F_q[[T]]\cdot X_0$. Proposition 2.14 implies that if $P\in I(X_0)$ then the set of zeroes of $P$ contains the set of zeroes of $P_1$ which is of dimension 2, hence is strictly greater than $\g X$. 
\medskip
Also, we see that the answer to Question 2.1 is "No". 
\medskip
Moreover, we see that any t-holonomic sequence $X=\sum_{i=0}^\infty x_iT^i$ has an invariant: the dimension of $V(I(X))$, and the set of $(i,j)$ such that the coefficient at $\tau^iT^j$ of $min(I(X))$ is non-zero. 
\medskip
{\bf Example.} Let $P=(\tau + a) + \tau^2T$, where $a\in \p$ is a parameter. It has one (up to multiplication by $\n F_q[[T]]$) solution; we denote it by $X=\sum_{i=0}^\infty x_iT^i$. 
\medskip
{\bf Lemma.} $V(X)=AP$.
\medskip
{\bf Proof.} To write. 
\medskip
To write explicitly equations satisfied by $x_i$, we start from the notion of an ordered partition $\pi$ of a number $n$. We have: $\pi=(\al_1, \dots, \al_k)$ such that $\al_i\in \n Z^{>0}$, and $\sum \al_i=n$. We attach to $\pi$ a polynomial $M_\pi:=\sum_{i=1}^k \al_i q^i+1$, and 
\medskip
{\bf Lemma.} $x_n^{q^2}=\sum_{i=0}^n$ (of what? to continue.)
\medskip
On $V(I(X))$ (to write).
\medskip
Let $X\in \p[[T]]$ be a holonomic sequence. Let us describe $V(I(X))$. We have: $I(X)$ is a saturated left ideal of $A$. Let $P_0:=min(I(X))$. We have:
$$V(P_0)=V(I(X)).$$
Hence, it is a $\n F_q[[T]]$-vector space, containing $X$, its dimension can be $>1$. We denote it by $E(X)$ (from "extension"). 
\medskip
For any $Y\in E(X)$ we have $E(Y)\subset E(X)$.
\medskip
{\bf Questions.} 1. Does exist a characterization of $X$ such that $\forall \ Y\in E(X) \ (Y\ne 0)$ we have $E(Y)=E(X)$? 
\medskip
2. Is it possible that $\exists \ X_1, \ X_2$ such that $E(X_1)\cap E(X_2)\ne 0$ and $E(X_1)\not\subset  E(X_2)$, $E(X_2)\not\subset  E(X_1)$?
\medskip
{\bf 3. Types of the stable formulas for $v(x_k)$
 for t-holonomic sequences.}
\medskip
For the classical holonomic sequences $x_0, \ x_1,\dots$
(D-case) it is well-known what
is the asymptotic behaviour of $x_k$ for $k\to\infty$, see for
example, [K], Theorem 2,
page 10. Here we prove propositions on the asymptotic behaviour of
$v(x_k)$ for t-
holonomic sequences for some particular cases, we formulate conjectures
for more
general cases, and state the general research problems.
\medskip
Let us prove the first proposition on the asymptotic behaviour of
$v(x_k)$. Let
$$P=\sum_{\ga=c_0}^{r_0}a_\ga\tau^\ga+\sum_{\ga=c_1}^{\vk_1}b_{1\ga}\tau^\ga T\in\p[T,\tau]\eqno{(3.1)}$$
be a polynomial having $n=1$, where $a_{r_0}, \ a_{c_0}, \ b_{1,\vk_1},\ b_{1,c_1}\ne0.$
\medskip
{\bf Proposition 3.2.} Let $x_0, \ x_1,\dots$ be a
t-holonomic sequence for $P$ satisfying
(3.1) such that $v(x_k)\to+\infty$. Then $c_1\ge c_0$, and $\exists \ \al>0, \ \be$ such that for all
sufficiently large $k$ we have:
$$\hbox{If $c_1>c_0$ then $v(x_k)=\al q^{(c_1-c_0)k}+\be$; \ \ if $c_1=c_0$ then $v(x_k)=\al k+\be$}\eqno{(3.3)}$$

{\bf Proof.} The $k$-th equation of the form (1.6.*) is
$$a_{r_0}x_k^{q^{r_0}}+a_{r_0-1}x_k^{q^{r_0-1}}+...+a_{c_0}x_k^{q^{c_0}}+t_k=0\eqno{(3.4)}$$
where $t_k$ (the tail) is
$$b_{1,\vk_1}x_{k-1}^{q^{\vk_1}}+b_{1,\vk_1-1}x_{k-1}^{q^{\vk_1-1}}+...+b_{1,c_1}x_{k-1}^{q^{c_1}}.$$
Here $x_{k-1}$ is already known, and $x_k$ is unknown. The Newton
polygon of (3.4) is the
lower convex hull of points $(q^j, \ v(a_j))$ for $j=c_0, \dots, r_0$, and the point $(0, v(t_k))$.
\medskip
Since $v(x_k)\to +\infty$, we have $v(t_k)\to +\infty$ as well, i.e. for all sufficiently large $k$
the segment $[(0, v(t_k)), (q^{c_0}, \ v(a_{c_0}))]$ belongs to the Newton polygon of (3.4). 
\medskip
The slopes of the lower convex hull of points $(q^j, \ v(a_j))$
for $j=c_0, \dots, r_0$ do not depend on $k$.
Let $C$ be the maximal value of $v(x)$, where $x$
is a root of (3.4) corresponding to these slopes. Since $v(x_k)\to +\infty$, 
for all sufficiently large
$k$ we have $v(x_k)>C$, i.e. for
these $k$ the root $x$ of (3.4) which is equal to $x_k$,
corresponds to the segment $[ (0, v(t_k)), (q^{c_0}, \ v(a_{c_0}))]$ of the Newton polygon of (3.4).
\medskip
Again, since $v(t_k)\to +\infty$, for all sufficiently
large $k$ we have $v(b_{1,c_1}x_{k-1}^{q^{c_1}})< v(b_{1,\ga}x_{k-1}^{q^{\ga}})$ for all $\ga>c_1$, i.e.
$$v(t_k)=v(b_{1,c_1}x_{k-1}^{q^{c_1}})=v(b_{1,c_1})+q^{c_1}v(x_{k-1}).$$

Since $v(x_k)$ corresponds to the segment $[(0, v(t_k)), (q^{c_0}, \ v(a_{c_0}))]$ of the
Newton polygon of (3.4), we get:
$$v(x_k)=\frac{v(b_{1,c_1})+q^{c_1}v(x_{k-1})-v(a_{c_0})}{q^{c_0}},$$
i.e. for all sufficiently large $k$ we have
$$v(x_k)=q^{c_1-c_0}v(x_{k-1})+const\eqno{(3.5)}$$
Obviously (3.5) implies (3.3). $\square$
\medskip
Recall ([GL21]) that elements $X$ such that $v(t_k)\to +\infty$ are called small solutions. They are used in calculation of $h_1,\ h^1$ of t-motives.
\medskip
{\bf Remark 3.6.} If $n$ of $P$ is $>$ 1
then it can happen that for all sufficiently large $k$ there will
be two or more terms in the tail of (1.6.*) having the same $v$.
For example, if $n=2$, $a_0, \ b_{10}, \ b_{20}\ne0$ and $v(b_{10})=v(b_{20})-const$, where $const$ is from (3.5), then for all
sufficiently large $k$ we have $v(b_{10}x_{k-1})=v(b_{20}x_{k-2})$. In this case, in order to find
$v(x_k)$, we need to consider the residues of terms of
(1.6.*) by the maximal ideal.
\medskip
Hence, we can state
\medskip
{\bf Open question 3.7.} Let $x_0, \ x_1, \dots$ be a
t-holonomic sequence such that $v(x_k)\to +\infty$. Is
it true that for all sufficiently large $k$ we have:
\medskip
$v(x_k)$ is a linear combination of $k^\al q^{\be k},\ \al, \be=0, \ 1, \dots$?
\medskip
We have an analog of Proposition 3.2 for growing t-holonomic sequences:
\medskip
{\bf Proposition 3.8.} Let $x_0, \ x_1, \dots$ be a
t-holonomic sequence for $P$ of the form (3.1) such that $v(x_k)\to -\infty$.
Then
$\vk_1\ge r_0$, and $\exists \ \al>0, \ \be$ such that for all
sufficiently large $k$ we have:
$$\hbox{If $\vk_1>r_0$ then $v(x_k)=-\al q^{(\vk_1-r_0)k}+\be$; \ \ if $\vk_1=r_0$ then $v(x_k)=-\al k+\be$}\eqno{(3.9)}$$

{\bf Proof} is similar to the proof of Proposition 3.2, we use its
notations. For all sufficiently large $k$ the Newton polygon of
(3.4) is one segment $[(0, v(t_k)); (q^{r_0}, v(a_{r_0}))]$. For all sufficiently large $k$ we have $v(t_k)=q^{\vk_1}v(x_{k-1})+v(b_{1,\vk_1})$. Hence, for all sufficiently large
$k$ we have:
$$v(x_k)=\frac{v(b_{1,\vk_1})+q^{\vk_1}v(x_{k-1})-v(a_{r_0})}{q^{r_0}},$$
i.e. for all sufficiently large $k$ we have
$$v(x_k)=q^{\vk_1-r_0}v(x_{k-1})+const\eqno{(3.10)}$$
Obviously (3.10) implies (3.9). $\square$
\medskip
Analog of 3.6 holds for this case. Also, we can state an analog of 3.7.
\medskip
Further, let $P\in \p[T,\tau]$, and let $\p[T,\tau]\cdot P$ be the left ideal generated by $P$. We denote $V(\p[T,\tau]\cdot P)$
simply by $V(P)$, it is the set of $X \in \p[[T]]$ such that $PX=0$. \ $V(P)$ is, as we wrote above, a free $\n F_q[[T]]$-module
of rank $r_0$.
\medskip
{\bf Open question 3.11.} Is it true that there exists a $\n F_q[[T]]$-basis of $V(P)$ such that for
all elements $X$ of this basis we have:
\medskip
A. Either $v(x_k)\to +\infty$;
\medskip
B. Either $v(x_k)\to$ a constant $\ne \pm\infty$;
\medskip
C. Or $v(x_k)\to -\infty$?
\medskip
Finally, we can ask what is an analog of the formula 3.7 for $X$ of
type B.
\medskip
We can also ask what is an analog of 3.11 for any left ideal of $\p[T,\tau]$
(not necessarily a principal ideal).
\medskip
{\bf 4. Properties of t-holonomic sequences.}
\medskip
Recall that we refer the ''classical'' holonomic sequences as the
D-case. We expect that the
analogs of closure properties of holonomic sequences in the D-case (see, for example, [K], Theorem 4; [KP], Theorem 7.2) hold in the t-case. Let us recall the definition and some properties for the D-case. Let $D: \ \n C[[T]]\to \n C[[T]]$ be the operator of derivation. 
\medskip
{\bf Definition 4.1.} A sequence $x_0, \ x_1, \dots$, where $x_i\in \n C$, is called a holonomic sequence, if two equivalent conditions hold: 
\medskip
{\bf 4.2.} $\exists \ r_1>0$ and polynomials $p_0(x),\dots,p_{r_1}(x)\ne0$ such that $\forall \ k\ge0$ we have
$$p_0(k)x_k+p_1(k)x_{k+1}+...+p_{r_1}(k)x_{k+r_1}=0$$

{\bf 4.3.} $\exists \ r_2 >0$ and polynomials $q_0(T),\dots,q_{r_2}(T)\ne0$ such that $f(T):=\sum_{i=0}^\infty x_iT^i\in \n C[[T]]$ satisfies $P(f)=0$, where $P=\sum_{i=0}^{r_2}q_iD^i$ is a differential operator, i.e. $$q_0f+q_1f'+...+q_{r_2}f^{(r_2)}=0\eqno{(4.3.1)}$$

Let $n_1$, resp. $n_2$ be the maximal degree of $p_i(x)$, resp. $q_i(T)$. We have ([KP], Th. 7.1):
\medskip
If (4.2) holds for $x_0, \ x_1, \dots$ with some $r_1$, $n_1$, then (4.3) holds for $x_0, \ x_1, \dots$ with $r_2\le n_1$, $n_2\le r_1+n_1$;
\medskip
If (4.3) holds for $x_0, \ x_1, \dots$ with some $r_2$, $n_2$, then (4.2) holds for $x_0, \ x_1, \dots$ with $r_1\le r_2+ n_2$, $n_1\le r_2$. 
\medskip
Let us consider the t-case. We see that $r_0$, resp. $n$ from (1.9) are analogs of $r_2$, resp. $n_2$ of (4.3). 
\medskip
What are analogs of $n_1$, resp. $r_1$ of (4.2)?
\medskip
For example, let $P$ be fixed. For the
t-case the set of $X$ satisfying (1.6) is a free $\n F_q[[T]]$-module of dimension $r_0$, and for the
D-case the set of $f(T)$ satisfying (4.3.1) is a free
$\n C$-module of dimension $r_2$ (for generic $P$).
\medskip
Nevertheless, there is no complete analogy between the t-case and the
D-case. For example, for the D-case there is no analog of the following
fact that holds in the t-case. Let $\g C$ be the Carlitz module and $x_0=0, \ x_1, \ x_2, \dots$ a sequence of elements
satisfying $x_{i+1}^q+\th x_{i+1}=x_i$. It is a t-holonomic
sequence. An analog of the $Exp$ map for the Carlitz module is the
ordinary $exp: \ \n C^+/2\pi i\to \n C^*$. Hence, an analog
of the above $\{x_n\}$ is a sequence $a_n:=exp(2\pi i/p^n)$. It is not holonomic.
\medskip
{\bf 4.4.} Let us mention here a research problem. For the D-case we
can define holonomic functions in several variables, see, for example,
[K], Definition 3:
\medskip
{\bf Definition 4.4.1.} Let $F(T) = \sum_Ix_IT^I$ where $T=(T_1,\dots,T_m)$ is a multivariable and $I \in \{0,1,\dots\}^m$ is a multiindex.
$F$ is called holomorphic if $\forall \ i\in \{1,\dots,m\}$
there exist non-zero polynomials $q_{ij}\in \n C[T]\ (j\in \{0,\dots,r_2\})$ such that
$$q_{i,r_2}\frac{\partial^{r_2}F}{\partial T_i^{r_2}}+ q_{i,r_2-1}\frac{\partial^{r_2-1}F}{\partial T_i^{r_2-1}}+...+q_{i1}\frac{\partial F}{\partial T_i}+q_{i0}F=0\eqno{(4.4.1.1)}$$

{\bf Open problem 4.4.2.} Does exist an analog of this definition for
the t-case?
\medskip
Since the analog of the differential operator is the Frobenius
automorphism, we need $m$ Frobenius automorphisms. We denote them
by $\tau_i, \ i=  1,\dots,m$. Hence, the first step to
solve 4.4.2 is
\medskip
{\bf Open problem 4.4.3.} Does exist a field (we denote it by
$\n C_{m,\infty}$) having $m$ commuting automorphisms $\tau_i$, such
that it is algebraically closed, complete with respect to all $m$
valuations related to $\tau_*$, and such that $\prod \tau_i$ is the Frobenius
$x\mapsto x^q$ (hence the subfield stable with respect to all $\tau_i$ is
$\n F_q$)?
\medskip
Let us try to construct $\n C_{m,\infty}$. We should consider not one element
$\th$ but elements $\th_1, \dots, \th_m$. A field $\n R_{m,\infty}$ is
maybe the quotient field of $\n F_q[[\th_1, \dots, \th_m]]$, or maybe there exists $\Omega$ satisfying $\n Z_{\ge0}^m\subset \Omega \subset \n Z^m$ such that $\n R_{m,\infty}$ is the set of
the sums $$\sum_{I\in\Omega}c_I\th^I$$
where $I$ is a multiindex, $\th:=\{\th_1, \dots, \th_m\}$, and $c_I\in\n F_q$.
\medskip
We have $m$ Frobenius endomorphisms $\tau_i, \ i=  1,\dots,m$ on $\n R_{m,\infty}$, namely $\tau_i(\th_i)=\th_i^q,\  \tau_i(\th_j)=\th_j$ for $i\ne j$.
\medskip
Now let us consider an algebraic closure of $\n R_{m,\infty}$ denoted
by $\overline{\n R_{m,\infty}}$. Now we should define $\n C_{m,\infty}$ which should be a completion of $\overline{\n R_{m,\infty}}$ with respect to all $m$ valuations.
\medskip
The analog of (4.4.1.1) for the t-case should be the following:
\medskip
Let $X=\sum_{i=0}^\infty x_iT^i\in \n C_{m,\infty}[[T]]$. $X$ is called t-m-holonomic, if $\forall \ i=1,\dots,m$ \ $\exists \ P_{ij}\in \n C_{m,\infty}[T] \ \ (j\le r)$ such that $$P_{ir}\tau_i^r(X)+P_{i,r-1}\tau_i^{r-1}(X)+...+P_{i1}\tau_i(X)+
P_{i0}X=0.\eqno{(4.4.4)}$$

{\bf Remark 4.4.5.} The above $P_{ij}$ must satisfy some relations.
Really, let us consider the simplest case $m=2,\ r=1, \ P_{11}=P_{21}=1$. Condition $\tau_1\circ\tau_2(x)=x^q$ and compatibility of (4.4.4) for $i=1, \ 2$ imply

$$\tau_2(P_{10})P_{20}=\tau_1(P_{20})P_{10}\eqno{(4.4.5)}$$

The analogous relation holds for the D-case. Really, if in (4.4.1.1) we
have $m=2,\ r_2=1, \ q_{11}=q_{21}=1$, then (4.4.1.1) implies
$$\frac{\partial q_{10}}{\partial T_2}=\frac{\partial q_{20}}{\partial T_1}$$

which is an analog of (4.4.5).
\medskip
Let us prove the closure property for the addition, for the t-holonomic
sequences:
\medskip
{\bf Theorem 4.5.} Let $x_0, x_1,\dots$ and $y_0, y_1,\dots$ be t-holonomic sequences. Then their memberwise sum $x_0+y_0, x_1+y_1,\dots$ is also a t-holonomic sequence.
\medskip
{\bf Proof.} Like in [K], [KP], we use a method of indeterminate coefficients, with some minor modifications. For the reader's convenience, we do not give the general proof, but we consider only some particular cases: the general idea is the same. First, we consider the simplest non-trivial case $r_0=n=1$, $\vk_1=0$ (notations of (1.5)) for both $x_0, x_1,\dots$, $y_0, y_1,\dots$. For their sum $x_0+y_0, x_1+y_1,\dots$ we have $n=2$, $r_0=3$, i.e. there exist $P_2, \ P_1, \ P_0\in \p[\tau]$ such that $\forall \ i $ we have $$P_2(x_{i+2}+y_{i+2})+P_1(x_{i+1}+y_{i+1})+P_0(x_{i}+y_{i})=0\eqno{(4.6)}$$
and deg($P_2)=3$. Moreover, we can choose $P_1$ and $P_0$ of degrees 2 and 1 respectively. 
\medskip
Let us consider a vector space $V$ of dimension 18 over $\p$ whose basis elements are 
$$\matrix
x_{i+2}^{q^3}, \ x_{i+2}^{q^2}, \ x_{i+2}^{q}, \ x_{i+2}, x_{i+1}^{q^2}, \ x_{i+1}^{q}, \ x_{i+1}, x_{i}^{q}, \ x_{i}, \\
 y_{i+2}^{q^3}, \ y_{i+2}^{q^2}, \ y_{i+2}^{q}, \ y_{i+2}, y_{i+1}^{q^2}, \ y_{i+1}^{q}, \ y_{i+1}, y_{i}^{q}, \ y_{i}\endmatrix\eqno{(4.7)}$$

(considered as abstract symbols). Conditions that $x_0, x_1,\dots$ and $y_0, y_1,\dots$ are holonomic sequences give us relations (here $a_*=a_*(x)$, $b_*=b_*(x)$ are from (1.5) for $x_0, \ x_1,\dots)$: 
$$a_1x_{i+1}^{q}+a_0x_{i+1}+b_{10}x_i=0$$
$$a_1^qx_{i+1}^{q^2}+a_0^qx_{i+1}^q+b_{10}^qx_i^q=0$$
$$a_1x_{i+2}^{q}+a_0x_{i+2}+b_{10}x_{i+1}=0\eqno{(4.8)}$$
$$a_1^qx_{i+2}^{q^2}+a_0^qx_{i+2}^q+b_{10}^qx_{i+1}^q=0$$
$$a_1^{q^2}x_{i+2}^{q^3}+a_0^{q^2}x_{i+2}^{q^2}+b_{10}^{q^2}x_{i+1}^{q^2}=0$$

and analogous relations for $y_*$, with coefficients $a_*(y)^{q^*}$, $b_*(y)^{q^*}$. There are 10 equalities. We can consider the left hand sides of these equalities as elements of $V$ (linear combinations of basis elements). The matrix of coefficients of these elements is of size $10\times18$ of rank 10, because $a_1(x), \ a_1(y)\ne0$, the matrix consists of two $5\times9$ diagonal blocks, both of them have a $5\times5$ triangilar submatrix with the diagonal entries $a_1(*)^{q^*}$. We denote by $V_0$ the subspace of $V$ generated by these elements.
\medskip
Let $V_1$ be the subspace of $V$ generated by vectors $x_{i+\al}^{q^\be}+y_{i+\al}^{q^\be}$ for all admissible combinations of $\al, \ \be$ from (4.7). We have dim $V_1=9$, dim $V_0=10$, hence $V_0\cap V_1\ne0$. 
\medskip
A non-zero element $v\in V_0\cap V_1\ne0$ gives us $P_2, \ P_1, \ P_0$ from (4.6). Really, let 
$$v=\sum_{\al, \be}c_{\al,\be}(x_{i+\al}^{q^\be}+y_{i+\al}^{q^\be})$$ (the same set of $\al, \ \be$). Then $P_\al:=\sum_\be c_{\al,\be}\tau^\be$ satisfy (4.6). 
\medskip
\medskip
Now, we consider the case when both $x_0, x_1,\dots$ and $y_0, y_1,\dots$ have $n=1$, the same $r$, and their $\vk_1$ is $r-1$ (the maximal possible value). 
\medskip
We have: $n$ of $\{x_i+y_i\}$ is 2. Let in (4.6) we have deg $P_2=\la$ (an unknown), deg $P_1=\la-1$, deg $P_0=\la-2$. 
\medskip
$V$ has a basis 
$$x_{i+2}^{q^\la}, \dots , x_{i+2}, \ \ x_{i+1}^{q^{\la-1}}, \dots, x_{i+1}, \ \ x_{i}^{q^{\la-2}}, \dots,  x_{i}\eqno{(4.9)}$$

and the same elements for $y_{i+*}^{q^{*}}$. Hence, we have: dim $V=6\la$, dim $V_1=3\la$. 
\medskip
Equations (4.8) have the form 

$$a_rx_{i+1}^{q^r}+...+a_0x_{i+1}+b_{1,r-1}x_i^{q^{r-1}}+...+b_{10}x_{i}=0$$
$$\dots$$
$$a_r^{q^{\la-r-1}}x_{i+1}^{q^{\la-1}}+...+a_0^{q^{\la-r-1}}x_{i+1}^{q^{\la-r-1}}+ b_{1,r-1}^{q^{\la-r-1}}x_{i}^{q^{\la-2}}+...+b_{10}^{q^{\la-r-1}}x_{i}^{q^{\la-r-1}}=0$$

($\la-r$ equations)
\medskip
and 

$$a_rx_{i+2}^{q^r}+... +a_0x_{i+2}+b_{1,r-1}x_{i+1}^{q^{r-1}} +...+    b_{10}x_{i+1}=0$$
$$\dots$$
$$a_r^{q^{\la-r}}x_{i+2}^{q^\la}+...+a_0^{q^{\la-r}}x_{i+2}^{q^{\la-r}}+b_{1,r-1}^{q^{\la-r}}x_{i+1}^{q^{\la-1}}+...+b_{10}^{q^{\la-r}}x_{i+1}^{q^{\la-r}}=0$$

($\la-r+1$ equations). Hence, we have: dim $V_0=4\la-4r+2$ and the condition $V_0\cap V_1\ne0$ is satisfied if $7\la-4r+2>6\la$. This means that for any $\la\ge4r-1$ there is a non-trivial linear dependence relation between $(x_{i+\al}+y_{i+\al})^{q^{\be_\al}}$, $\al =0,\ 1, \ 2$, $\be_0\le\la-2$, \ $\be_1\le\la-1$, \ $\be_2\le\la$. 
\medskip
The proof for the general case follows the same ideas. One more case is treated below. $\square$
\medskip
Let $x_0, x_1, \dots$ be a holonomic sequence. The set of $P\in \p[T,\tau]$ such that $P(\sum_{i=0}^\infty a_iT^i)=0$ is an ideal in $\p[T,\tau]$. By analogy with the D-case we can conjecture that it can be non-principal (compare with [GL21], Proposition 2.3.2 affirming that for a $T$-divisible submodule of $\p[[T]]$ the defining equation is unique). Let us consider an evidence to this fact. (To compare with the end of Section 2).
\medskip
{\bf 4.10.} Let $x_0, x_1, \dots$, resp. $y_0, y_1, \dots$ be two t-holonomic sequences satisfied by equations having $n=2$, $r_0=r$ arbitrary, $\vk_1, \vk_2 \le r-1$. Let us find $n$ and $r_0$ for their sum $x_0+y_0, \ x_1+y_1, \dots$. Acting as above we choose first $n=4$, $\la$ an indeterminate. The analog of (4.9) is
$$x_{i+4}^{q^\la}, \dots , x_{i+4}, \ \ x_{i+3}^{q^{\la-1}}, \dots, x_{i+3}, \ \ x_{i+2}^{q^{\la-1}}, \dots,  x_{i+2}, \ \ x_{i+1}^{q^{\la-2}}, \dots, x_{i+1}, \ \ x_{i}^{q^{\la-2}}, \dots,  x_{i}$$
and the same elements for $y_{i+*}^{q^{*}}$. Hence, we have: dim $V=10\la-2$, dim $V_1=5\la-1$. 
\medskip
Equations (4.8) have the form 

$$a_rx_{i+2}^{q^r}+...+a_0x_{i+2}+b_{1,r-1}x_{i+1}^{q^{r-1}}+...+b_{10}x_{i+1}+b_{2,r-1}x_{i}^{q^{r-1}}+...+b_{20}x_{i}=0$$
$$\dots\eqno{(4.11)}$$
$$a_r^{q^{\la-r-1}}x_{i+2}^{q^{\la-1}}+...+a_0^{q^{\la-r-1}}x_{i+2}^{q^{\la-r-1}}+ b_{1,r-1}^{q^{\la-r-1}}x_{i+1}^{q^{\la-2}}+...+b_{10}^{q^{\la-r-1}}x_{i+1}^{q^{\la-r-1}}+$$
$$b_{2,r-1}^{q^{\la-r-1}}x_{i}^{q^{\la-2}}+...+b_{20}^{q^{\la-r-1}}x_{i}^{q^{\la-r-1}}=0$$

($\la-r$ equations);

$$a_rx_{i+3}^{q^r}+...+a_0x_{i+3}+b_{1,r-1}x_{i+2}^{q^{r-1}}+...+b_{10}x_{i+2}+b_{2,r-1}x_{i+1}^{q^{r-1}}+...+b_{20}x_{i+1}=0$$
$$\dots\eqno{(4.12)}$$
$$a_r^{q^{\la-r-1}}x_{i+3}^{q^{\la-1}}+...+a_0^{q^{\la-r-1}}x_{i+3}^{q^{\la-r-1}}+ b_{1,r-1}^{q^{\la-r-1}}x_{i+2}^{q^{\la-2}}+...+b_{10}^{q^{\la-r-1}}x_{i+2}^{q^{\la-r-1}}+$$
$$b_{2,r-1}^{q^{\la-r-1}}x_{i+1}^{q^{\la-2}}+...+b_{20}^{q^{\la-r-1}}x_{i+1}^{q^{\la-r-1}}=0$$

($\la-r$ equations) and 
$$a_rx_{i+4}^{q^r}+...+a_0x_{i+4}+b_{1,r-1}x_{i+3}^{q^{r-1}}+...+b_{10}x_{i+3}+b_{2,r-1}x_{i+2}^{q^{r-1}}+...+b_{20}x_{i+2}=0$$
$$\dots\eqno{(4.13)}$$
$$a_r^{q^{\la-r}}x_{i+4}^{q^{\la}}+...+a_0^{q^{\la-r}}x_{i+4}^{q^{\la-r}}+ b_{1,r-1}^{q^{\la-r}}x_{i+3}^{q^{\la-1}}+...+b_{10}^{q^{\la-r}}x_{i+3}^{q^{\la-r}}+$$
$$b_{2,r-1}^{q^{\la-r}}x_{i+2}^{q^{\la-1}}+...+b_{20}^{q^{\la-r}}x_{i+2}^{q^{\la-r}}=0$$
($\la-r+1$ equations). Hence, we have: dim $V_0=6\la-6r+2$ and the condition $V_0\cap V_1\ne0$ is satisfied if $11\la-6r+1>10\la-2 \ \iff \ \la> 6r-3$. 
\medskip
Now we repeat the above arguments for $n=5$. 
\medskip
The analog of (4.9) is
$$x_{i+5}^{q^\la}, \dots , x_{i+5}, \ \ x_{i+4}^{q^{\la-1}}, \dots, x_{i+4}, \ \ x_{i+3}^{q^{\la-1}}, \dots,  x_{i+3}, \ \ x_{i+2}^{q^{\la-2}}, \dots, x_{i+2}, $$ $$x_{i+1}^{q^{\la-2}}, \dots, x_{i+1}, \ \ x_{i}^{q^{\la-3}}, \dots,  x_{i}$$
and the same elements for $y_{i+*}^{q^{*}}$. Hence, we have: dim $V=12\la-6$, dim $V_1=6\la-3$. 
\medskip
Analogs of (4.11) --- (4.13) are 4 sets of equations, with $\la-r-1$, $\la-r$, $\la-r$, $\la-r+1$ equations in these sets. Hence, we have: dim $V_0=8\la-8r$ and the condition $V_0\cap V_1\ne0$ is satisfied if $14\la-8r-3>12\la-6 \ \iff \ \la> 4r-\frac32$.  
\medskip
We see that the ideal of elements $P\in\p[T,\tau]$ such that $P[\sum_{i=0}^\infty (x_i+y_i)T^i]=0$ contains two elements: the first one has $n=4, \ r_0=6r-2$, the second one has $n=5, \ r_0=4r-1$. This is an evidence that this ideal is not principal: compare with [KP], Fig. 7.1, p. 143. 
\medskip
{\bf 4.14. }We can consider the following question. Let $x_0, \ x_1, \dots$ be a t-holonomic sequence such that all $x_i\ne0$, and $X:=\underset{i=0}\to{\overset{\infty}\to{\sum}}x_iT^i$. Obviously $\exists! \ P_0\in\p[[T]][\tau]$ of the form $$P_0=(\tau+b_0)+b_1T+b_2T^2+..., \ \ \ b_i\in\p\eqno{(4.15)}$$ such that $P_0X=0$. Is it true that this $P_0\in\p[T,\tau]$, i.e. is it true that if $\exists \ P\in\p[T,\tau]$ such that $PX=0$ then $\exists \ P_0\in\p[T,\tau]$ of the form (4.15) such that $P_0X=0$? 
\medskip
The answer is: not. We have a simple counter-example. Let $$P=\tau^2+a\tau+\th+\tau T$$ where $v(a)=-q$. The condition $PX=0$ is 
$$x_0^{q^2}+ax_0^q+\th x_0=0$$
$$x_1^{q^2}+ax_1^q+\th x_1+x_0^q=0$$
$$x_2^{q^2}+ax_2^q+\th x_2+x_1^q=0$$
etc. We find $x_i$ as functions of $a$; we have $v(x_i)=q^i+q^{i-1}+...+q+1$. 
\medskip
Further, the condition $P_0X=0$ is (here $b_i$ are unknowns)
$$x_0^q+b_0 x_0=0$$
$$x_1^q+b_0 x_1+b_1x_0=0$$
$$x_2^q+b_0 x_2+b_1x_1+b_2x_0=0$$
etc. By induction, we get that $v(b_i)=(i+1)q-1$, i.e. $P_0\in \p[[T]][\tau]-\p[T,\tau]$. 
\medskip
{\bf Remark.} Does exist $P_0\in \p[T,\tau]$ for $P$ coming from the equation of the Tate module $T_T(M)$, $M$ of rank 2, i.e. $P=\tau^2+a\tau+\th-T$? In this case the finding of $v(b_i)$ is more complicated, because of equality of $v$ of some terms. 
\medskip
{\bf 5. Systems of affine equations.}
\medskip
A system of $\la$ affine equations in $\mu$ variables $X_1, \dots, X_\mu\in \p[[T]]$ has the $i$-th equation ($i=1,\dots,\la$) of the form (see (0.3))
$$\sum_{j=1}^\mu P_{ij}X_j=0\eqno{(5.1)}$$
where coefficients $P_{ij}\in \p[T,\tau]$ are as above. They form a $\la\times \mu$-matrix denoted by $\g P$. 
\medskip
{\bf 5.1A.} Let $X_1, \dots, X_\mu\in \p[[T]]$ be fixed. Let $\g Q\in (\p[T,\tau])^m, \ \g Q=\{Q_1, \dots, Q_m\}$ where $Q_i \in \p[T,\tau]$ be such that $$\sum_{j=1}^\mu Q_iX_i=0.\eqno{(5.1A.1)}$$ 
The set of these $\g Q$ is a left $\p[T,\tau]$-submodule of $(\p[T,\tau])^m$; we denote it by $I(X_1,\dots,X_\mu)$. Study of its invariants is a research problem. We shall show in Proposition 5.8 that $I(X_1,\dots,X_\mu)$ contains elements with only one non-zero entry. 
\medskip
{\bf Remark 5.2.} We can consider a version of (5.1) for the case $P_{ij}\in \p[[T]]\{\tau\}$, as well as for $P_{ij}\in \p[[T]]\{\{\tau\}\}$.
\medskip
{\bf Remark 5.3.} For the D-case we have: let $a_{ij\al}$ be functions, 
$$D_{ij}:=\sum_\al a_{ij\al}\left(\frac{d}{d x}\right)^\al$$ differential operators, $X_j$ unknown functions. The $i$-th equation of an analog of (5.1) is $$\sum_{j=1}^\mu D_{ij}(X_j)=0.$$
\medskip
The set of solutions to (5.1) is a $\n F_q[[T]]$-submodule of $(\p[[T]])^\mu$. We denote it by $\g W=\g W(\g P)$. It is $T$-divisible, i.e. if $X=(X_1, \dots, X_\mu)\in\g W$ and $X/T\in (\p[[T]])^\mu$ then $X/T\in \g W$. See [GL21], (2.2.3) for the case $\la=\mu=1$. 
\medskip
Let $\la=\mu$. The $\n F_q[[T]]$-rank of $\g W$ for the simplest case is the following. We let $$P_{ij}=\sum_{\be, \ga\ge0}a_{ij;\be\ga}\tau^\ga T^\be$$ We denote by $r_{ij0}$ the maximal value of $\ga$ such that $a_{ij;0\ga}\ne0$. Further, we denote $X_i=\sum_{\be=0}^\infty x_{i\be}T^\be$. The $i$-th equation of (5.1) for the $T$-free terms is
$$R_i:=\sum_{j=1}^\mu\sum_{\ga=0}^{r_{ij0}}a_{ij;0\ga}\ x_{j0}^{q^\ga}=0$$

{\bf 5.4.} It is an additive polynomial in $x_{10},\dots, x_{\mu0}$ of degree $q^{m_i}$ where $m_i:=\underset{j}\to{\max} \ r_{ij0}$. Hence, if the intersections of the hypersurfaces $R_i=0$ are transversal at all points then the quantity of these intersection points is $q^{\sum_i m_i}$, and the $\n F_q[[T]]$-rank of $\g W$ is $\sum_{i=1}^\mu m_i$. 
\medskip
{\bf 5.5.} We denote $\Cal A:=(a_{ij;00})$. Since $\frac{\partial R_i}{\partial x_{j0}}=a_{ij;00}$ we get that if $|\Cal A|\ne0$ then the rank of $\g W$ is $\sum_i m_i$.
\medskip
Let $|\Cal A|=0$, i.e. the corank of $\Cal A$ is $\de>0$. We consider linear transformations of lines of $\Cal A$ making its first $\de$ lines the zero lines. We make the linear transformations of equations $R_i=0$ with the same coefficients. As a result, we get that the first $\de$ equations of the equations $R_i=0$  (i.e. the values $i=1,\dots, \de$) will be without free terms, hence there exist $\tilde R_i\in \p\{\tau\}$ such that $R_i=\tilde R_i^q$. The system 
$$\tilde R_1=0, \dots, \tilde R_{\de}=0, R_{\de+1}=0, \dots, R_{\mu}=0\eqno{(5.6)}$$ has the same set of solutions $\g W$, and its $m_i$ will be $$\tilde m_1:=m_1-1, \dots, \tilde m_{\de}:=m_{\de}-1, \ m_{\de+1}, \dots, m_{\mu}$$
and the problem of finding the rank of $\g W$ will be reduced to the same problem for the system (5.6). If necessary we can repeat the same constriction several times. Hence, the case of general $\g P$ can always be reduced to the case of $\g P$ having $|\Cal A|\ne0$. 
\medskip
{\bf Question 5.7.} Let $\g W\subset (\p[[T]])^\mu$ be a $T$-divisible $\n F_q[[T]]$-submodule. How to restitute $\g P$ by $\g W$? For the case $\la=\mu=1$ we can restitute $P$ uniquely (under some conditions), see [GL21], Proposition 2.3.2.

\medskip
Let now $\la=\mu=n$. For a generic $\g P$ we can eliminate $X_2,\dots, X_n$ using cofactors, like in the commutative case. 
\medskip
{\bf Proposition 5.8.} For a generic quadratic $\g P$ there exist non-zero cofactors $C_1,\dots, C_n\in \p[T,\tau]$ such that 
$$\forall \ j=2,\dots, n\hbox{ we have }\sum_{i=1}^n C_iP_{ij}=0.$$
\medskip
{\bf Proof.} It is analogous to the proof for the D-case, see 5.37. Let $k$ be the maximal degree of all $P_{ij}$ as polynomials in $\tau$. We let $$C_i=\sum_{j=0}^{k(n-1)} z_{ij}\tau^j\eqno{(5.9)}$$ where $z_{ij}\in \p[T]$ are indeterminate coefficients. Products $C_iP_{ij}$ are of degree $kn$ in $\tau$, hence for a fixed $j$ the equation $\sum_{i=1}^n C_iP_{ij}=0$ is a system of $kn+1$ linear equations in $z_{ij}$ with coefficients in $\p[T]$. Since  $j=2,\dots, n$, we have $(kn+1)(n-1)$ linear equations. There are $n(k(n-1)+1)$ unknowns $z_{ij}$ - one more than equations. 

Hence, we get a matrix denoted by $\Cal M=\Cal M(\g P)$ of size $$(kn+1)(n-1)\times n(k(n-1)+1)$$ with entries in $\p[T]$, of coefficients of this system of linear equations. The system is $\Cal M\cdot (z_{10},\dots,z_{n,k(n-1)})^t=0$ ($t$ means transposition). By definition, $z_{ij}$ are cofactors of $\Cal M$ and $C_i$ are from (5.9). $\square$
\medskip
{\bf Definition 5.10.} $\sum_{i=1}^n C_iP_{i1}$ is called the determinant of $\g P$ along the first column. Notation: $det_{1,c}(\g P)$. The same determinant for the $i$-th column is denoted by $det_{i,c}(\g P)$.
\medskip
{\bf Remark. } If we choose $k$ greater than the maximal degree of $P_{ij}$ as polynomials in $\tau$, then the last line of $\Cal M$ and hence all $C_i$ would be 0. 
\medskip
{\bf 5.11.} Clearly Proposition 5.8 and Definition 5.10 are valid for any commutative ring $R$ endowed with an automorphism $x\mapsto x^{(1)}$, and $R\{\tau\}$ satisfying $\tau x = x^{(1)}\tau$ for any $x\in R$. Namely, in notations of (5.1), $P_{ij}\in R\{\tau\}$, $X_i\in R$. In our case $R$ is $\p[T]$, and $X_i$ from (5.1) are from $\p[[T]]$. 
\medskip
We denote the projection of $\p[[T]]^n\to\p[[T]]$ to its first coordinate by $\al=\al_1$, and the set of $X\in \p[[T]]$ satisfying $$det_{1,c}(\g P)X=0$$ by $\g W_1=\g W_1(\g P)$. Obviously $\al(\g W)\subset \g W_1$.
\medskip
{\bf 5.12.} Under some reasonable conditions the map $\al|_\g W: \g W \to \g W_1$ is an isomorphism. Really, let $n=2$ and $P_{21}$,  $P_{22}$ are "coprime from the right", i.e. for $X\in \p[[T]]$ $$P_{21}X=P_{22}X=0\ \implies \ X=0.$$  Then obviously $\al|_\g W$ is a monomorphism. Let us consider epimorphisity. First, a $T$-divisible $\n F_q[[T]]$-submodule of $\p[[T]]$ of finite rank has no proper $T$-divisible submodules of the same rank, because its elementary divisors are powers of $T$, which contradicts to the condition of $T$-divisibility. Further, let all $m_i$ from (5.4) be equal to $k$ from 5.8, and for $\Cal A$ from (5.5) we have $|\Cal A|\ne0$. In this case the rank of $\g W$ is $kn$. Proof of 5.8 shows that for a generic case the rank of $\g W_1$ is also $kn$. It remains to prove that $\al(\g W)$ is $T$-divisible. We restrict ourselves by the case $n=2$. Let $(XT,Y)\in \g W$. We have $P_{12}\cdot Y, \ P_{22}\cdot Y\in \p[[T]]\cdot T$. For a generic case, $\sum_\ga a_{12;0\ga}\tau^\ga$ and $\sum_\ga a_{22;0\ga}\tau^\ga$ are coprime from the right. This implies that $Y\in \p[[T]]\cdot T$, hence $(X,Y/T)\in \g W$ and $X\in \al(\g W)$.
\medskip
Explicitly, we can describe $\al^{-1}$ as follows. Let us consider an equation $$PX=K\eqno{(5.16)}$$ where $P\in \p[T,\tau]$ is a coefficient, $X\in \p[[T]]$ an unknown and $K\in \p[[T]]$ a free term. Let $P=\sum_i P_iT^i$ where $P_i\in \p\{\tau\}$. If $P_0\ne0$ then (5.16) has solutions. 
\medskip

Let $n=2$ and $X_1\in \g W_1(\g P)$. To find $X_2$ such that $(X_1, \ X_2)\in \g W$, we should solve a system $$\matrix P_{12}X_2=-P_{11}X_1\\ P_{22}X_2=-P_{21}X_1\endmatrix\eqno{(5.17)}$$ Both these equations are of type (5.16). Since $X_1\in \g W_1$ and $\al|_\g W$ is an epimorphism, we get that they have a common root. 
\medskip
If $n>2$ then the system (5.17) has the form $$\matrix P_{12}X_2+...+P_{1n}X_n=-P_{11}X_1\\ \dots \\P_{n2}X_2+...+P_{nn}X_n=-P_{n1}X_1\endmatrix\eqno{(5.18)}$$ Let us fix $i\in [1,\dots,n]$. We denote by $\hat \g P_{i1}$ the submatrix of $\g P$ obtained by elimination of the $i$-th line and first column. We disregard the $i$-th equation of (5.18) and find $X_2$ satisfying the remaining equations using the method of elimination of $X_3, \dots,X_n$ of Proposition 5.8. Namely, we get an equation  $$det_{2,c}(\hat\g P_{i1})X_2=K\eqno{(5.19)}$$ where $K\in \p[[T]]$ is a linear combination of $P_{j1}X_1$ ($j\ne i$) with cofactors. Since $X_1\in \g W_1$ and $\al|_\g W$ (we hope) is an epimorphism, we get that all (5.19) have a common root. 
\medskip
Continuing the process we get $(X_1, X_2,\dots, X_n)\in \g W$ over $X_1$. 
\medskip
In some particular cases (see Section 9, $n=2$: for this case $P_{12}$ is invertible in a larger ring) it is easy to get a proof that in these cases $\al|_\g W$ is an isomorphism. A proof for a general case is a problem of further research. 
\medskip
{\bf 5.20. Scaling.} Let $f\in \p[[T]]$ be an element satisfying a relation $$C_1f^{(1)}=C_0f\eqno{(5.21)}$$ where $C_1, \ C_0\in \p[T]$ are some non-0 coefficients. Let $\g P, \ \g X$ be from (0.3). It is easy to modify $\g P$ in order to get an analog of (0.3) satisfied by $f\g X$, i.e. to find $\g P_m$ (subscript $m$ means "modified" such that $$\g P_m\cdot f\g X=0.\eqno{(5.22)}$$

Really, let, like in (5.8), $k$ be the maximal degree of all $P_{ij}$ as polynomials in $\tau$. For any $l=0,\dots,k$ we let $$C_{[0]}:=C_0^{(k-1)}\cdot C_0^{(k-2)}\cdot ... \cdot C_0;$$
$$C_{[l]}:=C_0^{(k-1)}\cdot C_0^{(k-2)}\cdot ... \cdot C_0^{(l)}\cdot C_1^{(l-1)}\cdot ... \cdot C_1^{(1)}\cdot C_1\ (l\ne0,\ k); \eqno{(5.23)}$$
$$C_{[k]}:=C_1^{(k-1)}\cdot C_1^{(k-2)}\cdot ... \cdot C_1.$$
Obviously (5.21) implies (because the ring $\p[T]$ is commutative) that $\forall \ l=0,\dots,k$ $$\hbox{ the elements $C_{[l]}\tau^l(f)$ are equal. }\eqno{(5.24)}$$

Let $\forall \ i, \ j \in [1,\dots, r]$ we denote $P_{ij}=\underset{l=0}\to{\overset{k}\to{\sum}}P_{ijl}\tau^l$, where $P_{ijl}\in \p[T]$. The formula for $(\g P_m)_{ij}$ --- the $(i,j)$-th entry of $\g P_m$ --- is the following: 
$$(\g P_m)_{ij}=\underset{l=0}\to{\overset{k}\to{\sum}}\ P_{ijl}\ C_{[l]}\ \tau^l\eqno{(5.25)}$$
(5.24) implies that (5.22) holds for this $\g P_m$.
\medskip
{\bf Remark 5.26.} 1. This construction will be used for the case $C_1=T-\th$, $C_0=1$, i.e. $f=\Xi$, or $C_0=T-\th$, $C_1=1$, i.e. $f=\Xi^{-1}$.

2. We can generalize it as follows. The above scaling corresponds to a scalar map $\g X\mapsto f\ \g X$. We can consider a diagonal map. Namely, let $f_*:=(f_1,\dots, f_r)$ be a row matrix, $f_i\in \p[[T]]$. The map $\g X \mapsto f_*\cdot \g X$ corresponds to a diagonal matrix acting on $\p[[T]]^r$. In order to get existence of $\g P_m$ such that $$\g P_m\cdot (f_*\cdot\g X)=0\eqno{(5.27)}$$
$f_1,\dots, f_r$ must satisfy the condition: $\exists \ C_{10}, \ C_{11}, \dots, \ C_{r0}, \ C_{r1}\in \p[T]$ such that $$C_{11}f_1^{(1)}=C_{10}f_1=C_{21}f_2^{(1)}=C_{20}f_2=\dots =C_{r1}f_r^{(1)}=C_{r0}f_r.\eqno{(5.28)}$$
Since $C_{[l]}$ of (5.23) depend on $k$, we denote them by $C_{[l]}(k)$, and we denote $C_{[l]}(k)$ from the pair $C_{\al1}, \ C_{\al2}$, where $\al=1,\dots, r$, by $D_{\al,[l]}(k)$: $$D_{\al,[l]}(k):=C_{\al0}^{(k-1)}\cdot C_{\al0}^{(k-2)}\cdot ... \cdot C_{\al0}^{(l)}\cdot C_{\al1}^{(l-1)}\cdot ... \cdot C_{\al1}^{(1)}\cdot C_{\al1}$$ Particularly, $$D_{\al,[k-1]}(k-1)=C_{\al1}^{(k-2)}\cdot C_{\al1}^{(k-3)}\cdot ... \cdot C_{\al1}^{(1)}\cdot C_{\al1}.$$
Analog of (5.23) is $$C_{\al,[l]}=D_{1,[k-1]}(k-1) \cdot...\cdot D_{\al-1,[k-1]}(k-1)\cdot D_{\al,[l]}(k)\cdot D_{\al+1,[k-1]}(k-1)\cdot...\cdot D_{r,[k-1]}(k-1)$$
We have an analog of (5.24): $\forall \ \al=1,\dots,r$, $\forall \ l=0,\dots,k$ $$\hbox{ the elements $C_{\al,[l]}\tau^l(f_\al)$ are equal }\eqno{(5.29)}$$
Hence, an analog of (5.25) is the following: 
$$(\g P_m)_{ij}=\underset{l=0}\to{\overset{k}\to{\sum}}\ P_{ijl}\ C_{j,[l]}\ \tau^l\eqno{(5.30)}$$
The above formulas show that (5.27) holds for this $\g P_m$. 
\medskip
{\bf Remark 5.31.} Most likely this construction can be generalized for any linear transformation $\g X \mapsto \vf(\g X)$, $\vf\in GL_r(\p[T])$. We shall not write down analogs of (5.28) and (5.30) for this case, because we do not know whether these formulas have a practical aplication, or not.
\medskip
{\bf 5.32.} Let us give the explicit form of $\Cal M$. We denote $P_{ij}=\sum_{\al=0}^k a_{ij\al}\tau^\al$ where $a_{ij\al}\in \p[T]$ are coefficients. $\Cal M$ is a block matrix: $\Cal M=(\Cal M)_{\be\ga}$ where $\be=2,\dots,n$, $\ga=1,\dots,n$, and all $(\Cal M)_{\be\ga}$ are of size $(kn+1)\times (k(n-1)+1)$. The $(\de,\nu)$-th entry of $(\Cal M)_{\be\ga}$ is 
$$[ \ (\Cal M)_{\be\ga}\ ]_{\de\nu}=a_{\ga\be,\de-\nu}^{(\nu-1)}\hbox{ if $\de-\nu\in[0,\dots,k]$ and $[ \ (\Cal M)_{\be\ga}\ ]_{\de\nu}=0$ otherwise, }$$where for $A=\sum_{i=0}^\infty y_iT^i\in \p[T]$ we denote $A^{(k)}:=\sum_{i=0}^\infty y_i^{q^k}T^i$. For example, we have (case $n\ge4$): $(\Cal M)_{\be\ga}=$
$$\left(\matrix a_{\ga\be0}&0&0&\dots&0 &0&\dots&0 &0 
\\ a_{\ga\be1}&a_{\ga\be0}^{(1)}&0&\dots&0  &0&\dots&0 &0 
\\ a_{\ga\be2}&a_{\ga\be1}^{(1)}&a_{\ga\be0}^{(2)}&\dots&0&0&\dots&0  &0 
\\ \dots&\dots&\dots&\dots&\dots&\dots&\dots&\dots&\dots
\\ a_{\ga\be k}&a_{\ga\be,k-1}^{(1)}&a_{\ga\be,k-2}^{(2)}&\dots&a_{\ga\be0}^{(k)} &0&\dots&0&0
\\
\\ 0&a_{\ga\be k}^{(1)}&a_{\ga\be,k-1}^{(2)}&\dots&a_{\ga\be1}^{(k)} &a_{\ga\be0}^{(k+1)}&\dots&0&0
\\ \dots&\dots&\dots&\dots&\dots&\dots&\dots&\dots&\dots
\\0&0&0&\dots&0 &0&\dots&a_{\ga\be0}^{(k(n-1)-1)} &0
\\
\\0&0&0&\dots&0 &0&\dots&a_{\ga\be1}^{(k(n-1)-1)} &a_{\ga\be0}^{(k(n-1))} 
\\
\\0&0&0&\dots&0 &0&\dots&a_{\ga\be2}^{(k(n-1)-1)} &a_{\ga\be1}^{(k(n-1))} 
\\ \dots&\dots&\dots&\dots&\dots&\dots&\dots&\dots&\dots
\\0&0&0&\dots&0 &0&\dots&a_{\ga\be k}^{(k(n-1)-1)} &a_{\ga\be,k-1}^{(k(n-1))} 
\\
\\0&0&0&\dots&0 &0&\dots&0&a_{\ga\be k}^{(k(n-1))} 
\endmatrix \right)\eqno{(5.33)}$$

An analog of (5.33) for the D-case has slighly other form, see (5.39) below. 
\medskip
{\bf Remark 5.34.} Recall the definition of the $p$-resultant (see, for example, [E08]). Let $P_1=\sum_{i=0}^na_i\tau^i, \ P_2=\sum_{i=0}^mb_i\tau^i\in \p\{\tau\}$. Then (rough form; really, the structure of this matrix is the same as the structure of ordinary resultant: there are $m$ lines of $a_*$'s and $n$ lines of $b_*$'s): 
$$R_p(P_1,P_2):=\left|\matrix a_{n}&\dots &a_{0}&0&\dots&0 
\\ \dots&\dots&\dots&\dots&\dots&\dots
\\ 0&\dots&0&a_{n}^{q^{m-1}}&\dots &a_{0}^{q^{m-1}}
\\b_{m}&\dots &b_{0}&0&\dots&0 
\\ \dots&\dots&\dots&\dots&\dots&\dots
\\ 0&\dots&0&b_m^{q^{n-1}}&\dots &b_{0}^{q^{n-1}}
\endmatrix \right|$$
We see that the above blocks (5.33) are blocks of $R_p$. 
\medskip
{\bf Example 5.35:} $n=2, \ k=1$, $$\g P=\left(\matrix A_0+A_1\tau & B_0+B_1\tau \\ C_0+C_1\tau & D_0+D_1\tau \endmatrix \right)$$ The corresponding determinant is
$$det_{1,c}(\g P)=\left|\matrix A_0 & B_0 \\ C_0 & D_0 \endmatrix \right|\cdot \left|\matrix B_0 & B_1 \\ D_0 & D_1 \endmatrix \right|^{(1)}+\eqno{(5.36)}$$
$$+\left[\ \left|\matrix A_1 & B_0 \\ C_1 & D_0 \endmatrix \right|\cdot \left|\matrix B_0 & B_1 \\ D_0 & D_1 \endmatrix \right|^{(1)}+\left|\matrix A_0 & B_1 \\ C_0 & D_1 \endmatrix \right|^{(1)}\cdot \left|\matrix B_0 & B_1 \\ D_0 & D_1 \endmatrix \right|\ \right]\tau+\left|\matrix A_1 & B_1 \\ C_1 & D_1 \endmatrix \right|^{(1)}\cdot \left|\matrix B_0 & B_1 \\ D_0 & D_1 \endmatrix \right|\tau^2$$
For the case $n=2$, $k>1$, $\g P=\left(\matrix A & B \\ C & D \endmatrix \right)\in M_{2\times2}\p\{\tau\}$ the $\tau$-free term of $det_{1,c}(\g P)$ is $\left|\matrix A_0 & B_0 \\ C_0 & D_0 \endmatrix \right|\cdot R_p(B,D)^{(1)}$. 
\medskip
{\bf Remark 5.37.} Recall that the analog of 5.35 (i.e. $n=2, \ k=1$) for the D-case is the following. Let $a_*, \ b_*, \ c_*, \ d_*$ be functions which are coefficients, and $X_1, X_2$ unknowns. Analog of (5.1) is
$$\matrix a_0  X_1 +a_1  X'_1 + b_0  X_2 +b_1  X'_2=0 \\ \\
c_0  X_1 +c_1  X'_1 + d_0  X_2 +d_1  X'_2=0\endmatrix\eqno{(5.38)}$$

In order to eliminate $X_2$, we apply a differential operator $g_0 + g_1 \frac{d}{d x}$, resp. $h_0 + h_1 \frac{d}{d x}$ to the first, resp. second equation of (5.38), where $g_0, \ g_1, \ h_0, \ h_1$ are indeterminates. We want that after summation $X_2$ disappears. We have
$$\Cal M=\left(\matrix b_0&b'_0&d_0&d'_0\\ b_1&b_0+b'_1&d_1&d_0+d'_1\\ 0&b_1&0&d_1\endmatrix\right)\eqno{(5.39)}$$
and $g_*, \ h_*$ satisfy $\Cal M\cdot (g_0, \ g_1, \ h_0, \ h_1)^t=0$, hence $g_*, \ h_*$ are cofactors of $\Cal M$. 
\medskip
The second order operator $$(g_0 + g_1 \frac{d}{d x})(a_0 + a_1 \frac{d}{d x})+(h_0 + h_1 \frac{d}{d x})(c_0 + c_1 \frac{d}{d x})$$ annulating $X_1$, is the 1-determinant of the system.
\medskip
{\bf 5.40. Relations with other quasideterminants. } Let $A=\left(\matrix a_{11} & a_{12} \\ a_{21} & a_{22}  \endmatrix \right)$ be any matrix with entries in a non-commutative ring. Its quasideterminant $\g D(A)$ is $a_{11} - a_{12}  a_{22}^{-1} a_{21} $ (see, for example, [GGRW05]). Let $C_1$, $C_2$ be from 1.5  ?, i.e. $C_1a_{12} +C_2a_{22} =0$. This implies $a_{12}  a_{22}^{-1}=-C_1^{-1}C_2$ and $\g D(A)=a_{11}+C_1^{-1}C_2a_{21}$ hence 
$$det_{1,c}(A)=C_1\cdot\g D(A)$$
Similar formulas exist for $A$ of size $>2$.
\medskip
\medskip
{\bf 6. Systems of affine equations associated to Anderson t-motives.} 
\medskip
Let $M$ be an Anderson t-motive of rank $r$. We choose a basis $f_*:=(f_1,\dots,f_r)^t $ of $M$ over $\p[T]$. We denote by $Q\in M_{r\times r}(\p[T])$ the matrix of the action of $\tau$ on $f_*$, i.e. $\tau f_*=Q f_*$. We denote $M[[T]]:=M\underset{\p[T]}\to{\otimes}\p[[T]]$. $\tau$ acts on $M[[T]]$. We identify $M[[T]]$ with $\p[[T]]^r$ using the basis $f_*$, and under this identification (here $Y$ is a row vector) $$M[[T]]^\tau=\{\ Y\in \p[[T]]^r \ |\  Y^{(1)}Q=Y\ \}\eqno{(6.1)}$$
(6.1) is a system of affine equations $$(Q^t\tau-I_r)Y^t=0\eqno{(6.2)}$$
According [A86], we have: $H^1(M)$ is the set of small solutions to (6.2) (recall that small solutions are solutions whose coefficients tend to 0, see Section 3 and [GL21], above Proposition 2.3.3. The small rank is the rank of the module of small solutions). 
\medskip
In order to solve (6.2), we use the above elimination process. Let $P=P_1(M):=det_{1,c}(Q^t\tau-I_r)$. We use notations of (1.5) for $P$. Conjecturally $r_0$ of (1.5) is $r$ --- the rank of $M$. 
\medskip
Also conjecturally, the small rank of solutions to (6.2) is equal to the small rank of solutions to $PX=0$. 
\medskip
Recall that $n$ from (1.5) is called the length of the tail of $P$. We conjectured in [GL21] that it is equal to the dimension of $M$. This is true for standard t-motives having $n=2$, $r=4$ or 5. We shall show below that this is true for Drinfeld modules and their 1-duals. But for elementary t-motives of dimension $\ge3$ this is wrong, see Section 9 below. 
\medskip
{\bf 6.3. $H_1(M)$}. All written above can be applied to $H_1(M)$ as well. Let 
$$\p\{T\}:=\{ \ \sum _{i=0}^\infty a_i T^i \in \p[[T]]\ |\ a_i\to0\ \}$$
We have $H_1(M)=\Hom_{\p[T]}(M,\p\{T\})^\tau$ which is the set of small elements in $\Hom_{\p[T]}(M,\p[[T]])^\tau$. 
We can identify column vectors $X\in \p[[T]]^r$ with elements of $\Hom_{\p[T]}(M,\p[[T]])$. Namely, let $\vf\in \Hom_{\p[T]}(M,\p[[T]])$. The column vector $X$ corresponding to $\vf$ is, by definition, $\left(\matrix \vf(f_1)\\ \dots\\ \vf(f_r) \endmatrix \right)$. The analog of (6.1) is 
$$\Hom_{\p[T]}(M,\p[[T]])^\tau=\{\ X\in \p[[T]]^r \ |\  QX=X^{(1)}\ \}\eqno{(6.3.1)}$$
(6.3.1) is a system of affine equations (analog of (6.2) ): $$(I_r\tau-Q)X=0\eqno{(6.3.2)}$$
\medskip
Let us recall a definition of $T_\g T(M)$ --- the $\g T$-Tate module of $M$. Let $E(M)$ be the t-module corresponding to $M$. If we fix a basis $e_1, \dots,e_n$ of $M$ over $\p\{\tau\}$ then we get an identification of $E(M)$ and $\p^n$. By definition, we have
$$T_\g T(M):=\{\ x_0=0, \ x_1, \ x_2, ... \in E(M)\ | \ \forall\ i \ T(x_{i+1})=x_i\  \}\eqno{(6.3.3)}$$

{\bf Proposition 6.3.4.} The system (6.3.2) essentially coincides with the system (6.3.3) defining the Tate module $T_\g T(M)$ (via equations $T(x_{i+1})=x_i$ ). 
\medskip
{\bf Remark 6.3.5.} To get an explicit form of equations $T(x_{i+1})=x_i$ defining $T_\g T(M)$ we use a matrix of multiplication by $T$ in $e_1, \dots,e_n$, while the above (6.3.2) uses a basis $f_1, \dots,f_r$ of $M$ over $\p[T]$. The word "essentially" in the statement of 6.3.4 means that the systems 6.3.2, 6.3.3 coincide after some elementary transformations. It is easy to prove Proposition 6.3.4 for $M$ of any explicit form. 
\medskip
{\bf 6.4. Duality and multiplication by the Carlitz module.} Let $M$ be an Anderson t-motive. Its 0-dual $M^\vee$ is defined for example in [HJ], it is a more general object that a t-motive of [A86]: t-motives of [A86] are effective t-motives of [HJ]. By definition, 
$$Q(M^\vee)=(Q(M)^t)^{-1}$$ We see that \{(6.2) for $M$\} is equal to \{(6.3.2) for $M^\vee$\}. This is in concordance with the canonical isomorphism $H^1(M)=H_1(M^\vee)$. 
\medskip
Let $\g C$ be the Carlitz module. We have $Q(M\otimes\g C^{\otimes n})=(T-\th)^nQ(M)$. Let $\Xi=\sum_{i=0}^\infty\xi_iT^i$ be a solution to 
$$\Xi=(T-\th)\Xi^{(1)}$$satisfying $\forall \ i \ \ v(\xi_{i+1})>v(\xi_i)$
(see [G96], Example 5.9.36, p. 172). It is unique up to a factor from $\n F_q^*$. We have $\Xi, \ \Xi^{-1} \in \p\{T\}$.
\medskip
Formula (6.1) implies that we have an isomorphism $H^1(M)\to H^1(M\otimes\g C^{\otimes n})$ defined by $Y\mapsto \Xi^nY$. It is canonical, i.e. it does not depend on a choice of $f_1,\dots,f_r$. The same is true for $H_1(M)\to H_1(M\otimes\g C^{\otimes n})$. 
\medskip
Nevertheless, systems of affine equations (6.2) for $M$ and $M\otimes\g C^{\otimes n}$ are quite different. For example, let $M'$ be the 1-dual of $M$ (see [Tg95]; [GL07]). By definition, $M'=M^\vee\otimes\g C$; its advantage is that $M'$ is an effective t-motive if $M$ is pure and has the nilpotent operator $N=0$. Let us consider the case of a Drinfeld module $M$. We compare the affine equations for $H^1(M')$ (see (8.2), (8.4)) and the affine equation for $H_1(M)$, which coincides with the affine equation for $H^1(M^\vee)$. We see that they are essentially different. Compare also (9.1) and (9.3) --- case of $M$ of rank 4, dimension 2. Hence the operation of tensoring by $\g C$, although it leaves $H^1$ invariant, modifies the affine equation.  
\medskip
We can ask what will occur with affine equations after other tensor products and Hom's, not only after tensoring by $\g C$. 
\medskip
{\bf 6.6. Problems of further research. A.} For any types of t-motives we should find the set of bidegrees $\Cal D$, i.e. the set of pairs $(\be,\ga)\in (\n Z^+)^2$ such that $b_{\be\ga}\ne0$ (see (1.5)), as functions of entries of $Q$. 
\medskip 
{\bf B.} Is it possible to define duality on the set of systems of affine equations (or on its quotient set, by some appropriate equivalence relation) coming from the duality of Anderson t-motives, or, the same, duality between $H^1(M), \ H_1(M)$? See (7.5.2), (7.5.3) as an example of such duality for Drinfeld modules of rank 2. 
\medskip 
{\bf C.} We can find more analogs of the theorems of the D-case. For example, the member-by-member product of two holonomic sequences is also a holonomic sequence (D-case). Is the analogous property true for the case of t-motives? Namely, let $x_0, x_1,\dots$ and $y_0, y_1,\dots$ be solutions to equations $P_1X=0$, $P_2X=0$ (see Remark 1.1.1), where $P_1$, $P_2$ are of finite tail. Is the sequence $x_0y_0, x_1y_1,\dots$ a solution of such equation? 
\medskip
{\bf D. Galois action.} Let $K$ be a finite extension of $\n F_q(\th)$. If coefficients of $P_{ij}$ from (5.1) are in $K$ then the Galois group $\Gal(K):=\Gal(\bar K/K)$ acts on $\g W$ --- the set of solutions to (5.1). Hence, if $\la=\mu=r$ we have a map 
$$\rho: \ \Gal(K)\to GL_r(\n F_q[[T]])\eqno{(6.6.D.1)}$$
\medskip
{\bf Remark 6.6.D.2.} $\Gal(K)$ does not preserve small elements. Really, there are many inclusions $i: \bar K \hookrightarrow \p$. The restriction of the valuation $v$ on $\p$ with respect to an inclusion induces a valuation $i^*(v)$ on $\bar K$. 
\medskip
Let us fix an inclusion $i_0: \bar K \hookrightarrow \p$. An element $\sum_{j=0}^\infty a_jT^j$, where $a_j\in \bar K$, is small, if $\underset{j\to\infty}\to{\lim}v(i_0(a_j))=+\infty$. Equivalently, $\underset{j\to\infty}\to{\lim}[i_0^*(v)](a_j)=+\infty$. 
\medskip
For $\ga\in\Gal(K)$ we have $[\ga(i_0^*(v))](\ga(v))=[i_0^*(v)](a_j)$, hence $\underset{j\to\infty}\to{\lim}[\ga(i_0^*(v))](\ga(v))=+\infty$ but not necessarily $\underset{j\to\infty}\to{\lim}[i_0^*(v)](\ga(v))=+\infty$, i.e. 
$\Gal(K)$ does not preserve small elements. 
\medskip
We can conjecture that the standard properties of $\g p$-adic representations hold for $\rho$. Namely:
\medskip
{\bf Conjecture 6.6.D.3.} For almost all $P_{ij}$ we have: im $\rho$ is of finite index in $GL_r(\n F_q[[T]])$.
\medskip
Further, let $\g p$ be a prime of $K$. We denote by $\bar \g p$ an extension of $\g p$ to $\bar K$ and by $G_i\subset G_d\subset Gal(K)$ the inertia and decomposition subgroups. 
\medskip
Like for the case of the Galois action on Tate modules, we expect that we have:
\medskip
{\bf Conjecture 6.6.D.4.} For almost all $\gp$ \ \  $G_i$ acts trivially on $\g W$, hence the  action of $Fr_\g p$ on $\g W$ is defined up to a conjugation. 
\medskip
What is known on its characteristic polynomial? 
\medskip 
{\bf E.} Let $\g p\in \n F_q[T]$ be an irreducible element. Elements of the Tate module $T_\gp(M)$ are sequences $x_0, \ x_1,\dots$, where $x_i\in E(M)=\p^n$, $x_0=0$, such that $\g p(x_i)=x_{i-1}$. They are defined by systems of affine equations. We can define $T_\gp(M)_s$. What is their small rank? Does it depend on $\g p$? 
\medskip
As a first example, we can consider the case of $\g p$ of degree 1, i.e. $\g p=T+c$ for $c\in \n F_q$. 
\medskip
Further, the identification of $T_T(M)_s$ and $H_1(M)$ shows that there is a natural inclusion of $T_T(M)_s$ to $Lie(M)$. Does there exist something similar for $T_\gp(M)_s$?
\medskip
{\bf 7. Case of Drinfeld modules.} 
Let $M$ be a Drinfeld module of rank $r$. It is free of rank 1 over $\p\{\tau\}$. Let $e$ be the only element of a basis of $M$ over $\p\{\tau\}$ such that the action of $T$ on $e$ is given by the formula (here $a_i\in \p$ are parameters):
$$Te=\th e+a_1\ \tau e +a_2\ \tau^2 e +...+a_{r-1}\ \tau^{r-1} e + \tau^{r} e$$
We have: a column $f_*=(e,\ \tau e, \ \tau^2 e, \dots,\tau^{r-1} e)^t $. Its $Q$ is 
$$Q=\left(\matrix 0&1&0&0&\dots&0&0
\\ 0&0&1&0&\dots&0&0 
\\ \dots&\dots&\dots&\dots&\dots&\dots&\dots
\\  0&0&0&0&\dots&1&0
\\  0&0&0&0&\dots&0&1
\\ T-\th&-a_1&-a_2&-a_3&\dots&-a_{r-2}&-a_{r-1}\endmatrix \right)$$
and its $\g P=Q^t\tau-I_r$ is
$$\g P=\left(\matrix -1&0&0&\dots&0&0&-\th\tau+\tau T
\\ \tau&-1&0&\dots&0&0 &-a_1\tau
\\ 0&\tau&-1&\dots&0&0 &-a_2\tau
\\ \dots&\dots&\dots&\dots&\dots&\dots&\dots
\\  0&0&0&\dots&\tau&-1&-a_{r-2}\tau
\\  0&0&0&\dots&0&\tau&-1-a_{r-1}\tau\endmatrix \right)$$
We find $det_{r,c}(\g P)$ "manually", i.e. by explicit elimination of unknowns. We denote (compare with [GL21], Section 3): $Y=(y_0 \ y_1 \ \dots \ y_{r-1})$ (caution: shift of indices), where $y_i\in \p[[T]]$. (6.1) becomes 
$$\matrix (T-\th)y_{r-1}^{(1)}=y_0
\\ -a_1y_{r-1}^{(1)}+y_0^{(1)}=y_1
\\ -a_2y_{r-1}^{(1)}+y_1^{(1)}=y_2
\\ \dots
\\ -a_{r-1}y_{r-1}^{(1)}+y_{r-2}^{(1)}=y_{r-1}\endmatrix\eqno{(7.1)}$$
We denote (like in [GL21], the line below (3.8)) $$y_{r-1}=x_0+x_1T+x_2T^2+...,\eqno{(7.2)}$$ where $x_i\in \p$, and we substitute this formula to the equations of (7.1):
$$y_0=-\th x_0^q+(-\th x_1^q+x_0^q)T+(-\th x_2^q+x_1^q)T^2+(-\th x_3^q+x_2^q)T^3+...$$
$$y_1=(-\th^q x_0^{q^2}-a_1x_0^q) +(-\th^q x_1^{q^2}-a_1x_1^q+x_0^{q^2})T+(-\th^q x_2^{q^2}-a_1x_2^q+x_1^{q^2})T^2+$$
$$+(-\th^q x_3^{q^2}-a_1x_3^q+x_2^{q^2})T^3+...$$
$$y_2=(-\th^{q^2} x_0^{q^3}-a_1^qx_0^{q^2}-a_2x_0^{q}) +(-\th^{q^2} x_1^{q^3}-a_1^qx_1^{q^2}-a_2x_1^q+x_0^{q^3})T+$$ $$(-\th^{q^2}x_2^{q^3}-a_1^qx_2^{q^2}-a_2x_2^q+x_1^{q^3})T^2+(-\th^{q^2} x_3^{q^3}-a_1^qx_3^{q^2}-a_2x_3^q+x_2^{q^3})T^3+...$$ $$\dots$$
$$y_{r-1}=(-\th^{q^{r-1}} x_0^{q^r}-a_1^{q^{r-2}}x_0^{q^{r-1}}-a_2^{q^{r-3}}x_0^{q^{r-2}}-...-a_{r-2}^qx_0^{q^2}-a_{r-1}x_0^{q}) +$$
$$(-\th^{q^{r-1}} x_1^{q^r}-a_1^{q^{r-2}}x_1^{q^{r-1}}-a_2^{q^{r-3}}x_1^{q^{r-2}}-...-a_{r-2}^qx_1^{q^2}-a_{r-1}x_1^{q}\ + \ x_0^{q^r})T +\eqno{(7.3)}$$
$$(-\th^{q^{r-1}} x_2^{q^r}-a_1^{q^{r-2}}x_2^{q^{r-1}}-a_2^{q^{r-3}}x_2^{q^{r-2}}-...-a_{r-2}^qx_2^{q^2}-a_{r-1}x_2^{q}\ + \ x_1^{q^r})T^2 +...$$
Comparing (7.2) and (7.3) we get an affine equation to find $x_0, \ x_1, \ x_2,\dots$:
$$\th^{q^{r-1}} x_0^{q^r}+a_1^{q^{r-2}}x_0^{q^{r-1}}+a_2^{q^{r-3}}x_0^{q^{r-2}}+...+a_{r-2}^qx_0^{q^2}+a_{r-1}x_0^{q}+x_0=0\hbox{ (the head);}\eqno{(7.4)}$$
$$\th^{q^{r-1}} x_1^{q^r}+a_1^{q^{r-2}}x_1^{q^{r-1}}+a_2^{q^{r-3}}x_1^{q^{r-2}}+...+a_{r-2}^qx_1^{q^2}+a_{r-1}x_1^{q}+x_1\underset{\hbox{tail}}\to{\underbrace{-\ x_0^{q^r}}}=0;$$
$$\dots$$
and for any $i$
$$\th^{q^{r-1}} x_i^{q^r}+a_1^{q^{r-2}}x_i^{q^{r-1}}+a_2^{q^{r-3}}x_i^{q^{r-2}}+...+a_{r-2}^qx_i^{q^2}+a_{r-1}x_i^{q}+x_i\underset{\hbox{tail}}\to{\underbrace{-\ x_{i-1}^{q^r}}}=0\eqno{(7.4.1)}$$
hence $$det_{r,c}(\g P)=[\ \sum_{i=0}^r a_{r-i}^{q^{i-1}}\tau^i\ ] -\tau^rT$$ (here $a_0=\th$, $a_r=1$). 
We see that the length of the tail is 1; the tail consists of one term $-\ x_{i-1}^{q^r}$. The set of  bidegrees $\Cal D$  from 2.5.A is $(0,0);\  (0,1);\dots; (0,r);\ (1,r)$. 
\medskip
{\bf 7.5.} The affine equation for $H_1(M)$ coincides with the affine equation for $T_\g T(M)$ (Proposition 6.3.4). It is (we omit details of calculation): 
$$x_i^{q^r}+a_{r-1}x_i^{q^{r-1}}+...+a_{2}x_i^{q^{2}}+a_{1}x_i^{q}+\th x_i-x_{i-1}=0\eqno{(7.5.1)}$$
Apparently it differs essentially from the above equation (7.4.1). For example, for $r=2$, $i=0$ the equations (7.5.1), resp. (7.4.1) are
$$ x_0^{q^2}+a_1x_0^q+\th x_0=0\eqno{(7.5.2)}$$ and $$ \th^qx_0^{q^2}+a_1x_0^q+ x_0=0\eqno{(7.5.3)}$$
Equations (7.5.3) and (7.5.2) are different: there exist $a_1\in \n F_q(\th)$ such that (7.5.3) has a non-trivial root in $\n F_q(\th)$ while (7.5.2) has no such root, and vice versa. For example, there exist $a_1$ such that the Galois groups of (7.5.2), (7.5.3) are different. 
\medskip
{\bf 8. Case of duals of Drinfeld modules.} Let $M$ be as above. Its 1-dual $M'$ has rank $r$ and dimension $r-1$. We have $Q(M')=(T-\th)[Q(M)^{-1}]^t$:
$$Q(M')=\left(\matrix a_1& T-\th&0&0&\dots&0&0
\\ a_2&0& T-\th&0&\dots&0&0 
\\ \dots&\dots&\dots&\dots&\dots&\dots&\dots
\\  a_{r-2}&0&0&0&\dots& T-\th&0
\\  a_{r-1}&0&0&0&\dots&0& T-\th
\\ 1&0&0&0&\dots&0&0\endmatrix \right)$$ and 
$$\g P(M')=\left(\matrix -1+a_1\tau&a_2\tau&a_3\tau&\dots&a_{r-2}\tau&a_{r-1}\tau&\tau
\\ -\th\tau+\tau T&-1&0&\dots&0&0 &0
\\ 0&-\th\tau+\tau T&-1&\dots&0&0 &0
\\ \dots&\dots&\dots&\dots&\dots&\dots&\dots
\\  0&0&0&\dots&-\th\tau+\tau T&-1&0
\\  0&0&0&\dots&0&-\th\tau+\tau T&-1\endmatrix \right)$$
As above we make an elimination "manually". We denote (shift in 1 in indices): $Y=(y_1 \ y_2 \ \dots \ y_{r})$. (6.1) becomes 
$$a_1y_{1}^{(1)}+a_2y_{2}^{(1)}+...+a_{r-1}y_{r-1}^{(1)}+y_{r}^{(1)}=y_1$$
$$(T-\th)y_{1}^{(1)}=y_2$$
$$\dots\eqno{(8.1)}$$
$$(T-\th)y_{r-2}^{(1)}=y_{r-1}$$
$$(T-\th)y_{r-1}^{(1)}=y_r$$
We denote $N_k:=(T-\th^{q^{k}})(T-\th^{q^{k-1}})\cdot ... \cdot (T-\th^{q^2})(T-\th^q)$, $N_0:=1$, and we denote its coefficients by $v_{ki}$ (they are polynomials in $\th$): 
$$N_k:=\sum_{i=0}^kv_{ki}T^i$$
(8.1) implies $y_k^{(1)}=N_{k-1}\ y_1^{(k)}$.  We denote like in (3.2): $y_{1}=x_0+x_1T+x_2T^2+...$. The first equation of (8.1) becomes the following constituents of an affine equation: 
$$v_{r-1,0}x_0^{q^r}+a_{r-1}v_{r-2,0}x_0^{q^{r-1}}+a_{r-2}v_{r-3,0}x_0^{q^{r-2}}+...+a_{2}v_{10}x_0^{q^{2}}+a_{1}v_{00}x_0^{q}-x_0=0$$ (the head --- the 0-th equality of the affine equation); 
$$v_{r-1,0}x_1^{q^r}+a_{r-1}v_{r-2,0}x_1^{q^{r-1}}+a_{r-2}v_{r-3,0}x_1^{q^{r-2}}+...+a_{2}v_{10}x_1^{q^{2}}+a_{1}v_{00}x_1^{q}-x_1+$$
$$v_{r-1,1}x_0^{q^r}+a_{r-1}v_{r-2,1}x_0^{q^{r-1}}+a_{r-2}v_{r-3,1}x_0^{q^{r-2}}+...+a_3v_{21}x_0^{q^{3}}+a_2v_{11}x_0^{q^{2}}=0$$
(the first equality of the affine equation); 
$$v_{r-1,0}x_2^{q^r}+a_{r-1}v_{r-2,0}x_2^{q^{r-1}}+a_{r-2}v_{r-3,0}x_2^{q^{r-2}}+...+a_{2}v_{10}x_2^{q^{2}}+a_{1}v_{00}x_2^{q}-x_2+$$
$$v_{r-1,1}x_1^{q^r}+a_{r-1}v_{r-2,1}x_1^{q^{r-1}}+a_{r-2}v_{r-3,1}x_1^{q^{r-2}}+...+a_3v_{21}x_1^{q^{3}}+a_2v_{11}x_1^{q^{2}}+$$
$$v_{r-1,2}x_0^{q^r}+a_{r-1}v_{r-2,2}x_0^{q^{r-1}}+a_{r-2}v_{r-3,2}x_0^{q^{r-2}}+...+a_4v_{32}x_0^{q^{4}}+a_3v_{22}x_0^{q^{3}}=0$$
(the second equality of the affine equation); etc.
The first equation of the maximal length of the tail is the $(r-1)$-th equality in the affine equation: 
$$v_{r-1,0}x_{r-1}^{q^r}+a_{r-1}v_{r-2,0}x_{r-1}^{q^{r-1}}+a_{r-2}v_{r-3,0}x_{r-1}^{q^{r-2}}+...+a_{2}v_{10}x_{r-1}^{q^{2}}+a_{1}v_{00}x_{r-1}^{q}-x_{r-1}+$$
$$v_{r-1,1}x_{r-2}^{q^r}+a_{r-1}v_{r-2,1}x_{r-2}^{q^{r-1}}+a_{r-2}v_{r-3,1}x_{r-2}^{q^{r-2}}+...+a_3v_{21}x_{r-2}^{q^{3}}+a_2v_{11}x_{r-2}^{q^{2}}+$$
$$v_{r-1,2}x_{r-3}^{q^r}+a_{r-1}v_{r-2,2}x_{r-3}^{q^{r-1}}+a_{r-2}v_{r-3,2}x_{r-3}^{q^{r-2}}+...+a_4v_{32}x_{r-3}^{q^{4}}+a_3v_{22}x_{r-3}^{q^{3}}+$$ $$\dots $$ 
$$v_{r-1,r-3}x_{2}^{q^r}+a_{r-1}v_{r-2,r-3}x_{2}^{q^{r-1}}+a_{r-2}v_{r-3,r-3}x_{2}^{q^{r-2}}+$$
$$+v_{r-1,r-2}x_1^{q^r}+a_{r-1}v_{r-2,r-2}x_{1}^{q^{r-1}}+$$
$$+v_{r-1,r-1}x_0^{q^r}=0$$ Finally, we write the $k$-th equality in the affine equation for $k\ge r-1$: 
$$v_{r-1,0}\ x_{k}^{q^r}+a_{r-1}\ v_{r-2,0}\ x_{k}^{q^{r-1}}+a_{r-2}\ v_{r-3,0}\ x_{k}^{q^{r-2}}+...+a_{2}\ v_{10}\ x_{k}^{q^{2}}+a_{1}\ v_{00}\ x_{k}^{q}-x_{k}+$$
$$v_{r-1,1}\ x_{k-1}^{q^r}+a_{r-1}\ v_{r-2,1}\ x_{k-1}^{q^{r-1}}+a_{r-2}\ v_{r-3,1}\ x_{k-1}^{q^{r-2}}+...+a_3\ v_{21}\ x_{k-1}^{q^{3}}+a_2\ v_{11}\ x_{k-1}^{q^{2}}+$$
$$v_{r-1,2}\ x_{k-2}^{q^r}+a_{r-1}\ v_{r-2,2}\ x_{k-2}^{q^{r-1}}+a_{r-2}\ v_{r-3,2}\ x_{k-2}^{q^{r-2}}+...+a_4\ v_{32}\ x_{k-2}^{q^{4}}+a_3\ v_{22}\ x_{k-2}^{q^{3}}+$$ $$\dots \eqno{(8.2)}$$ 
$$+v_{r-1,r-3}\ x_{k-r+3}^{q^r}+a_{r-1}\ v_{r-2,r-3}\ x_{k-r+3}^{q^{r-1}}+a_{r-2}\ v_{r-3,r-3}\ x_{k-r+3}^{q^{r-2}}+$$
$$+v_{r-1,r-2}\ x_{k-r+2}^{q^r}+a_{r-1}\ v_{r-2,r-2}\ x_{k-r+2}^{q^{r-1}}+$$
$$+v_{r-1,r-1}\ x_{k-r+1}^{q^r}=0, \hbox{    or}$$
$$\sum_{i=0}^{r-1}\ \ \ \sum_{j=0}^{r-1-i}\ \ a_{r-j}\ v_{r-1-j,i}\ \ x_{k-i}^{q^{r-j}}=x_k$$
Here $a_r=1$, $i$ is the number of a line in (8.2), and $j$ the number of a term in the $i$-th line, both $i,\ j$ are counted from 0. 
Equivalently, we have $$det_{1,c}(\g P)=[\ \sum_{i=0}^{r-1}\ \ \ \sum_{j=0}^{r-1-i}\ \ a_{r-j}\ v_{r-1-j,i}\ \ \tau^{r-j}T^i \ ]-1$$
The set of bidegrees $\g D$ is the union of a triangle with vertices $(0,1); \ (0,r); \ (r-1,r)$ and a point (0,0). 
We see that the length of the tail is $r-1$. 
\medskip
{\bf 8.3.} Let us transform (8.2) to the form where $\vk_\be<r$ (here $\vk_\be$ are from (1.1)?), in order to show that the length of the tail of the corresponding affine equation is $r-1$. We restrict ourselves by the case $r=3$. Formula (8.2) for $k=0$ is (we have $v_{00}=v_{11}=v_{22}=1$):
$$v_{20}\ x_{0}^{q^3}+a_{2}\ v_{10}\ x_{0}^{q^2}+a_{1}\ x_{0}^{q}-x_{0}=0,\hbox{ i.e.}$$
$$x_{0}^{q^3}=-a_{2}\ v_{10}\ v_{20}^{-1}\ x_{0}^{q^2}-a_{1}\ v_{20}^{-1}\ x_{0}^{q}+v_{20}^{-1}\ x_{0}$$ Formula (4.2) for $k=1$ is 
$$v_{20}\ x_{1}^{q^3}+a_{2}\ v_{10}\ x_{1}^{q^2}+a_{1}\ x_{1}^{q}-x_{1}+$$
$$+v_{21}\ x_{0}^{q^3}+a_{2}\ x_{0}^{q^2}=0, \hbox{ or}$$

$$v_{20}\ x_{1}^{q^3}+a_{2}\ v_{10}\ x_{1}^{q^2}+a_{1}\ x_{1}^{q}-x_{1}+$$
$$+v_{21}\  ( -a_{2}\ v_{10}\ v_{20}^{-1}\ x_{0}^{q^2}-a_{1}\ v_{20}^{-1}\ x_{0}^{q}+v_{20}^{-1}\ x_{0})  +a_{2}\ x_{0}^{q^2}=0, \hbox{ or}$$

$$v_{20}\ x_{1}^{q^3}+a_{2}\ v_{10}\ x_{1}^{q^2}+a_{1}\ x_{1}^{q}-x_{1}+$$
$$+a_{2}(-\ v_{10}\ v_{20}^{-1}v_{21}+1)\ x_{0}^{q^2}-a_{1}\ v_{20}^{-1}v_{21}\ x_{0}^{q}+v_{20}^{-1}v_{21}\ x_{0}  =0, \hbox{ i.e.}$$

$$x_{1}^{q^3}=-v_{20}^{-1}[a_{2}\ v_{10}\ x_{1}^{q^2}+a_{1}\ x_{1}^{q}-x_{1}+$$
$$+a_{2}(-\ v_{10}\ v_{20}^{-1}v_{21}+1)\ x_{0}^{q^2}-a_{1}\ v_{20}^{-1}v_{21}\ x_{0}^{q}+v_{20}^{-1}v_{21}\ x_{0} ] $$ Formula (8.2) for $k=2$ is 
$$v_{20}\ x_{2}^{q^3}+a_{2}\ v_{10}\ x_{2}^{q^2}+a_{1}\ x_{2}^{q}-x_{2}+$$
$$+v_{21}\ x_{1}^{q^3}+a_{2}\ x_{1}^{q^2}+ x_{0}^{q^3}=0, \hbox{ or}$$

$$v_{20}\ x_{2}^{q^3}+a_{2}\ v_{10}\ x_{2}^{q^2}+a_{1}\ x_{2}^{q}-x_{2}+$$
$$-v_{21} v_{20}^{-1}[a_{2}\ v_{10}\ x_{1}^{q^2}+a_{1}\ x_{1}^{q}-x_{1}+$$
$$+a_{2}(-\ v_{10}\ v_{20}^{-1}v_{21}+1)\ x_{0}^{q^2}-a_{1}\ v_{20}^{-1}v_{21}\ x_{0}^{q}+v_{20}^{-1}v_{21}\ x_{0} ] $$ $$       +a_{2}\ x_{1}^{q^2}-a_{2}\ v_{10}\ v_{20}^{-1}\ x_{0}^{q^2}-a_{1}\ v_{20}^{-1}\ x_{0}^{q}+v_{20}^{-1}\ x_{0}       =0, \hbox{ or}$$

$$v_{20}\ x_{2}^{q^3}+a_{2}\ v_{10}\ x_{2}^{q^2}+a_{1}\ x_{2}^{q}-x_{2}+$$
$$a_{2}(-v_{21} v_{20}^{-1}\ v_{10}+1)\ x_{1}^{q^2}-a_{1}\ v_{21} v_{20}^{-1}\ x_{1}^{q}+v_{21} v_{20}^{-1}\ x_{1}+\eqno{(8.4)}$$
$$+a_{2}(\ v_{10}\ v_{20}^{-2}v_{21}^2-v_{21} v_{20}^{-1}\ -\ v_{10}\ v_{20}^{-1})\ x_{0}^{q^2}+a_{1}\ ( v_{20}^{-2}v_{21}^2-v_{20}^{-1})\ x_{0}^{q}+(v_{20}^{-2}v_{21}^2+v_{20}^{-1})\ x_{0} =0$$

Substituting values of $v_{10},\ v_{20},\ v_{21}$ we get that the coefficients of the last line of (8.4) are non-0, hence the corresponding affine equation has the length of tail 2. 
\medskip
{\bf 9. The elementary t-motives.} 
\medskip
Let us consider the following Anderson t-motives $M$ of dimension $n$, rank $r=2n$. Let $e_*=(e_1,\dots, e_n)^t$ be a basis of $M$ over $\p\{\tau\}$. The action of $T$ on $e_*$ is defined by the formula 
$$Te_*=\theta \ e_*+A\ \tau\ e_* +\tau^2\ e_*.$$ where $A\in M_{n\times n}(\p)$ is a matrix parameter. This t-motive is denoted by $M(A)$. For $n=1$ it is a Drinfeld module of rank 2. We shall call them elementary t-motives. 
\medskip
We choose a basis $f_*=(e_1,\dots, e_n, \tau e_1,\dots, \tau e_n)^t$ of $M$ over $\p[T]$. The matrix $Q$ in $f_*$ is $\left(\matrix 0&I_n\\ (T-\th)\ I_n &-A\endmatrix \right)$ (block form, entries are $n\times n$-matrices). Since, according (6.2), $\g P$ is $Q^t\tau-I_{2n}$, we have for this case $$\g P=\g P_{2n}=\left(\matrix -I_n& (T-\th)\tau\ I_n \\ \tau I_n & -A^t\tau-I_n\endmatrix \right)$$

Making elementary transformations (their formalization is an exercise) we reduce $\g P_{2n}$ to the following matrix $\g P_{1,n}:=((T-\th^q)\tau^2-1)\cdot I_n-A^t\tau$, i.e.

$$\g P_{1,n}=\left(\matrix (T-\th^q)\tau^2-a_{11}\tau-1&-a_{21}\tau&\dots&-a_{n1}\tau
\\ -a_{12}\tau&(T-\th^q)\tau^2-a_{22}\tau-1&\dots&-a_{n2}\tau
\\ \dots&\dots&\dots&\dots
\\ -a_{1n}\tau&-a_{2n}\tau&\dots&(T-\th^q)\tau^2-a_{nn}\tau-1\endmatrix \right)$$
Notations of Proposition 5.8 for this case are: $$C_1=\sum_{j=0}^{2n-2}y_{1j}\tau^j, \hbox { and for } i>1 \ \ \ C_i=\sum_{j=1}^{2n-3}y_{ij}\tau^j$$
The matrices $\Cal M$ are: 
$$\hbox{ for } n=2\ \ \  \Cal M_2=\left(\matrix -a_{21}&0&0&-1
\\ 0&-a_{21}^q&0&-a_{22}^q
\\ 0&0&-a_{21}^{q^2}&T-\th^{q^2}
\endmatrix \right)$$

and for $ n=3 \ \ \ \Cal M_3$ is the following matrix (it is drawn in two parts; its block structure is denoted by spaces):
$$\left(\matrix -a_{21}&0&0&0&0&&-1&0&0
\\ 0&-a_{21}^q&0&0&0&&-a_{22}^q&-1&0
\\ 0&0&-a_{21}^{q^2}&0&0&&T-\th^{q^2}&-a_{22}^{q^2}&-1
\\ 0&0&0&-a_{21}^{q^3}&0&&0&T-\th^{q^3}&-a_{22}^{q^3}
\\ 0&0&0&0&-a_{21}^{q^4}&&0&0&T-\th^{q^4}
\\
\\ -a_{31}&0&0&0&0&&0&0&0
\\ 0&-a_{31}^q&0&0&0&&-a_{32}^q&0&0
\\ 0&0&-a_{31}^{q^2}&0&0&&0&-a_{32}^{q^2}&0
\\ 0&0&0&-a_{31}^{q^3}&0&&0&0&-a_{32}^{q^3}
\\ 0&0&0&0&-a_{31}^{q^4}&&0&0&0&&\endmatrix \right)$$

$$\left(\matrix 0&0&0\\
-a_{23}^q&0&0\\
0&-a_{23}^{q^2}&0\\
0&0&-a_{23}^{q^3}\\
0&0&0\\ \\
-1&0&0\\
-a_{33}^q&-1&0\\
T-\th^{q^2}&-a_{33}^{q^2}&-1\\
0&T-\th^{q^3}&-a_{33}^{q^3}\\
0&0&T-\th^{q^4}\endmatrix \right)$$

For $n=2$ we have 

$$det_{1,c}(\g P_2)=\th^{q^3+q^2}\ a_{21}^{q+1}\ \tau^4
+(\th^{q^2}\ a_{11}^{q^2}\ a_{21}^{q+1}+\th^{q^2}\ a_{21}^{q^2+1}\ a_{22}^q)\ \tau^3+$$ 
$$+(a_{11}^q\ a_{21}^{q^2+1}\ a_{22}^q\ -\ a_{12}^q\ a_{21}^{q^2+q+1}\  +\ \th^{q^2}\ a_{21}^{q+1}\ +\ \th^{q}\ a_{21}^{q^2+q})\ \tau^2 \ +$$ $$+(a_{21}^{q^2+1}\ a_{22}^q\ +\ a_{11}\ a_{21}^{q^2+q}) \ \tau\ +\ a_{21}^{q^2+q}+\eqno{(9.1)}$$ $$+[\ -(\th^{q^3}+\th^{q^2})\ a_{21}^{q+1}\tau^4-(a_{11}^{q^2}\ a_{21}^{q+1}+a_{21}^{q^2+1}\ a_{22}^q)\tau^3-(a_{21}^{q^2+q}+a_{21}^{q+1})\tau^2\ ]\ T+$$ $$+a_{21}^{q+1}\ \tau^4 \ T^2$$
We see that this is [GL21], (3.9), up to a constant factor. (3.11) of [GL21] is the linear combination $C_1P_{11}+C_2P_{21}$.
\medskip
The first equation of $\g P_2\left(\matrix X_1\\ X_2\endmatrix \right)=0$ gives us $$X_2=[a_{21}^{-1}( (T-\th^q)\tau^2-a_{11}\tau-1)X_1]^{(-1)}.$$ This shows that for this case $\al|_\g W: \g W \to \g W_1$ (see (5.12) ) is an isomorphism. 
\medskip
{\bf Remark 9.2.} Analog of (9.1) for $H_1(M)$ is 
$$det_{1,c}=\ a_{12}^{q+1}\ \tau^4
+(a_{11}^{q^2}\ a_{12}^{q+1}+ a_{12}^{q^2+1}\ a_{22}^q)\ \tau^3+$$ 
$$+(a_{11}^q\ a_{12}^{q^2+1}\ a_{22}^q\ -\ a_{21}^q\ a_{12}^{q^2+q+1}\  +\ \th^{q^2}\ a_{12}^{q+1}\ +\ \th^{q}\ a_{12}^{q^2+q})\ \tau^2 \ +$$ 
$$+(\th^q\ a_{12}^{q^2+1}\ a_{22}^q\ +\ \th^q\ a_{11}\ a_{12}^{q^2+q}) \ \tau\ +\ \th^{q+1}\ a_{12}^{q^2+q}+\eqno{(9.3)}$$ $$+[\ -(a_{12}^{q^2+q}+ a_{12}^{q+1}) \ \tau^2-(a_{11}\ a_{12}^{q^2+q}+a_{12}^{q^2+1}\ a_{22}^q)\tau-(\th^q+\th)a_{12}^{q^2+q}\ ]\ T+$$ $$+a_{12}^{q^2+q}\ T^2$$

Description of a symmetry between (9.1) and (9.3) is an exercise. This is an analog for $n=2$ of a symmetry between (7.5.2) and (7.5.3). First, we should transpose $A$, because $M(A)'=M(A^t)$. 
\medskip
{\bf Research problems.} We described in [GL07] a class of t-motives called standard-2 t-motives. They are t-motives for which bases of $M$ over both $\p[T], \ \p\{\tau\}$ have a simple description. We want to find the set of bidegrees $\Cal D$ (see 6.6A) for their $det_{1,c}(\g P)$.
\medskip
Further, there is a conjecture that for almost all t-motives $M$ of rank $r$ we have $h^1(M)$ is either 0 or $r$. The above theory of elimination permits us to verify this conjecture by explicit calculations. 
\medskip

{\bf 10. Some examples.}
\medskip
Here we give some examples related to the Section 1A.
\medskip
{\bf Example 10.1.} Let $d = 2, \ P = (\tau+a)+(\tau+b)T$, $a, \ b\in \n F_{q^2}$, i.e. in notations 
of (1.5) we have $r_0=n=\vk_1 = 1$. Equations (1.6.*) for
this case become 

$$x_0^q+ax_0=0; \ \ \ \ \matrix x_1^q+ax_1+\\ x_0^q+bx_0\endmatrix =0; \ \ \ \ \matrix x_2^q+ax_2+\\ x_1^q+bx_1\endmatrix =0; \ \ \ \ \hbox{etc.}\eqno{(10.2)}$$

Let us show that $\exists \ C\in \n F_q(T)$ such that
$$Fr(x_0+x_1T+x_2T^2+...):=x_0^{q^2}+x_1^{q^2}T+x_2^{q^2}T^2+...=C(x_0+x_1T+x_2T^2+...)$$
(1A.3) implies (here $N, \ Tr$ are the norm and the trace from $\n F_{q^2}$ to $\n F_q$; all formulas for $n$ are proved immediately by induction):
$$x_0^q=-ax_0; \ \ x_1^q=-ax_1+(a-b)x_0; \ \ x_2^q=-ax_2+(a-b)x_1-(a-b)x_0;$$
$$ x_n^q=-ax_n+(a-b)x_{n-1}-(a-b)x_{n-2}+...+(-1)^n(a-b)x_0;$$
$$x_0^{q^2}=N(a)x_0; \ \ x_1^{q^2}=N(a)x_1+[-2N(a)+Tr(ab^q)]x_0;$$
$$x_2^{q^2}=N(a)x_2+[-2N(a)+Tr(ab^q)]x_1+[3N(a)-2Tr(ab^q)+N(b)]x_0;$$
$$x_3^{q^2}=N(a)x_3+[-2N(a)+Tr(ab^q)]x_2+[3N(a)-2Tr(ab^q)+N(b)]x_1+$$
$$[-4N(a)+3Tr(ab^q)-2N(b)]x_0;$$
$$\dots$$
$$x_n^{q^2}=N(a)x_n+\sum_{k=1}^n(-1)^k[(k+1)N(a)-k \cdot Tr(ab^q)+(k-1)N(b)]x_{n-k}$$ This means that 
$$C=N(a)+\sum_{k=1}^\infty(-1)^k[(k+1)N(a)-k \cdot Tr(ab^q)+(k-1)N(b)]T^k=$$ $$=\frac{N(a)+Tr(ab^q)T+N(b)T^2}{(T+1)^2}.$$

{\bf Example 10.3.} Let us give an explicit foirmula for the Lemma
1A.1.4. Let $d=3, \ r=2$ and $P=a_2\tau^2+a_1\tau+a_0$. Then (here $N$ is the norm from $\n F_{q^3}$ to $\n F_q$)
$$Q=a_2\cdot a_2^{q^2}\tau^4-a_1^q\cdot a_2^{q^2}\tau^3+(a_1\cdot a_1^{q}-a_0^qa_2)\tau^2-a_0^{q^2}a_1\tau+a_0^q\cdot a_0^{q^2};$$
$$QP=N(a_2)\tau^6+[N(a_1)-N(a_0a_1^qa_2^{q^2})]\tau^3+N(a_0).$$
Apparently $Q$ cannot be represented as a product of two
polynomials in $\tau$ of degree 2, with coefficients coming from
$a_i$ and cyclotomic units.
\medskip
\medskip
{\bf References}
\medskip
[A86] Anderson, G.W. t-motives. Duke Math. J., 1986, vol. 53, No. 2,
p. 457 --502.
\medskip
[E08] Eri\'c, Aleksandra Lj. The resultant of non-commutative
polynomials. Mat. Vesnik 60 (2008), no. 1, 38.
\medskip
[GGRW05] Gelfand, Israel; Gelfand, Sergei; Retakh, Vladimir; Wilson,
Robert Lee. Quasideterminants. Adv. Math. 193 (2005), no. 1, 56 -- 141.
\medskip
[G96] Goss, D. Basic structures of function field arithmetic.
Springer-Verlag, Berlin, 1996. xiv+422 pp.
\medskip
[GL07] Grishkov A., Logachev, D. Duality of Anderson t-motives.
2007. http://arxiv.org/pdf/math/0711.1928.pdf
\medskip
[GL20] Grishkov A., Logachev, D. Introduction to Anderson t-motives:
a survey. 2020. https://arxiv.org/pdf/2008.10657.pdf
\medskip
[GL21] Grishkov A., Logachev, $h^1\ne h_1$ for
Anderson t-motives. J. of Number Theory. 2021, vol. 225, p. 59 -- 89.
https://arxiv.org/pdf/1807.08675.pdf
\medskip
[HJ20] Hartl U.; Juschka A.-K. Pink's theory of Hodge structures and
the Hodge conjecture over function fields. ''t-motives: Hodge
structures, transcendence and other motivic aspects'', Editors G.
B\"ockle, D. Goss, U. Hartl, M. Papanikolas, European Mathematical
Society Congress Reports 2020.
\medskip
[K] Kauers, M. The Holonomic Toolkit.

https://www3.risc.jku.at/publications/download etc.
\medskip
[KP] Kauers, Manuel; Paule, Peter. The concrete tetrahedron. Texts
Monogr. Symbol. Comput. Springer Wien, New York, Vienna, 2011, x+203 pp.
\medskip
[NP] Namoijam Ch., Papanikolas M. Hyperderivatives of periods and
quasi-periods for Anderson t-modules.
https://arxiv.org/pdf/2103.05836.pdf
\medskip
[O33.1] Ore, O. On a special class of polynomials. Trans. Amer.
Math. Soc. 35 (1933), no. 3, 559 -- 584.
\medskip
[O33.2] Ore, O. Theory of non-commutative polynomials, Ann. of Math.
34 (1933), 480 -- 508.
\medskip
[Tg95] Taguchi, Yuichiro. A duality for finite t-modules. J. Math.
Sci. Univ. Tokyo 2 (1995), no. 3, 563 -- 588.
\enddocument